\documentclass[12pt,aps,pre,preprint,showkeys]{revtex4}

\usepackage{comment}
\usepackage{amsmath}
\usepackage{amssymb}
\usepackage{hyperref}
\usepackage[dvips]{graphicx}
\usepackage{colordvi}
\usepackage{times}
\usepackage{float}
\usepackage[a4paper,top=2.5cm,bottom=2.5cm,left=2cm,right=2cm,headsep=11pt]{geometry}
\newcommand{\be}{\begin{equation}}
\newcommand{\ee}{\end{equation}}
\newcommand{\bea}{\begin{eqnarray}}
\newcommand{\eea}{\end{eqnarray}}

\newcommand{\zl}{\lambda}

\newcommand {\bdm} {\begin{displaymath}}
\newcommand {\edm} {\end{displaymath}}
\newcommand {\ba}  {\begin{array}}
\newcommand {\ea}  {\end{array}}

\begin{document}

\title{Role of Ergodicity in the Transient Fluctuation Relation and a New Relation for a Dissipative Non-Chaotic Map}

\author{Paolo A. Adamo \email{adamatrice@gmail.com}}
\affiliation{Dipartimento di Scienze Matematiche, Politecnico di Torino, Corso Duca degli Abruzzi 24, 10129 Torino, Italy}

\author{Matteo Colangeli \email{colangeli@mat.ufmg.br}}
\affiliation{Departamento de Matem{\'a}tica-ICEx, UFMG, CP 702 Belo Horizonte - MG, 30161-970 Brazil}

\author{Lamberto Rondoni \email{lamberto.rondoni@polito.it}}
\affiliation{Dipartimento di Scienze Matematiche, Politecnico di Torino, Corso Duca degli Abruzzi 24, 10129 Torino, Italy.\\
INFN, Sezione di Torino, Via P. Giuria 1, 10125, Torino, Italy}

\keywords{Ergodicity; Equilibrium; Fluctuation Relations; Dissipation.}

\begin{abstract}

Deterministic dynamical systems such as the baker maps are useful to shed light on some of the  
conditions verified by deterministic models in non-equilibrium statistical physics. We investigate a $2$D  
dynamical system, enjoying a weak form of reversibility, with peculiar basins of attraction and steady states. 
In particular, we test the conditions required for the validity of the Transient Fluctuation Relation. 
Our analysis illustrates by means of concrete examples why ergodicity of the equilibrium dynamics seems to be a	 necessary condition for the Transient Fluctuation Relation to hold. This investigation then leads to the numerical verification 
of a kind of transient relation which, differently from the usual Transient Fluctuation Relation (FR), holds only asymptotically. 
At the same time, it is not a steady state fluctuation relation, because no fluctuations are present in the steady state.

\end{abstract}

\maketitle

%\pacs{...}

%\submitto{\JPA}

\section{Introduction}
The theory of Fluctuation Relations (FRs) originated with the seminal works of Evans, Cohen and Morris \cite{ECM} which were followed by many other works, beginning with \cite{ES} for the transient FR and \cite{GCa} for the steady state FR. These relations became increasingly popular in statistical mechanics, as they describe the statistical properties of systems even far from equilibrium \cite{MeRo07,BPRV,ESR2}. Therefore it is important to clarify the assumptions under which such relations can be derived.
Numerous works investigate the assumption that the invariant probability measure  is smooth along the unstable direction, and that the dynamics is time reversible \cite{CR,CKDR,PVTW,GRS}, in the derivation of the steady state FR for the phase space contraction rate $\Lambda$. However, the identification of the minimal requirements for the transient FR to hold is still an open question. To this purpose, we consider a map $L$ of the unit square $\mathcal{U}$ which depends on one parameter $\ell$ which can be tuned to produce either non-dissipative dynamics, in which case we speak of ``equilibrium'' dynamics, or dissipative dynamics, in which case we speak of ``non-equilibrium'' dynamics. Furthermore, $L$ verifies a weak notion of time reversibility, meaning that, given an initial condition $\overrightarrow{x}=(x,y) \in \mathcal{U} $ which corresponds to a given  phase space contraction $\Lambda_n$ in $n$ iterations, there exists an initial condition $\overrightarrow{x'}$ corresponding to the opposite phase space 
contraction, $-\Lambda_n$.
\begin{comment}This initial condition $\overrightarrow{x}$ is obtained from the application of a map $\tilde{G}$, which is not an involution on the square, as detailed below, hence it cannot be considered in the definition of time reversal.
\end{comment}

Our results can be summarized as follows:
\begin{itemize}
\item the transient FR is not verified by $L$, which is not time reversal invariant (although it is  ``weakly'' reversible) and is not ergodic at equilibrium. This illustrates why equilibrium ergodicity and time reversibility had to be invoked in the original derivations of the transient FR.
\item we nevertheless verify a kind of {\em transient asymptotic} FR which is not the steady state FR. A similar situation was envisaged in Ref.\cite {ESR2}, albeit for reversible dynamics. This seems to be a novel property for deterministic non-chaotic dynamics.
\end{itemize}

%%%%%%%%%%%%%%%%%%%%%%%%%%%%

\section{The map}
\label{sec:sec1}
Let us consider maps of the square, which, in spite of the lack of immediate physical application, 
are analytically tractable, hence allow a detailed analysis of the dynamics \cite{K05,RTV}.
Here, we introduce a generalized baker-like transformation, which corresponds to a modified version of the map previously discussed in Refs. \cite{CR,CKDR}.
 
\begin{figure}[htbp!]
\begin{center}	
\includegraphics[width=0.92\textwidth]{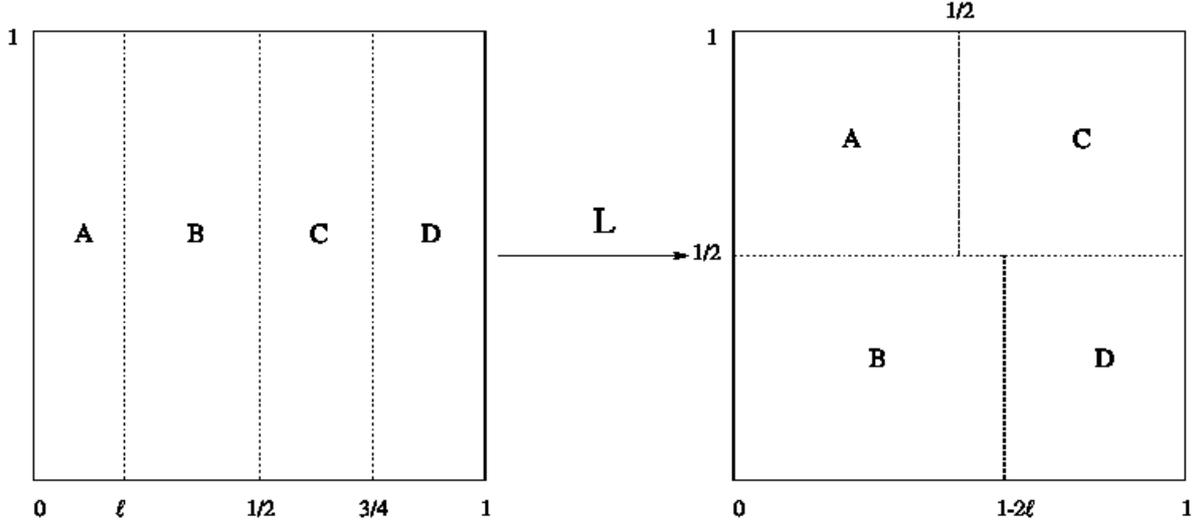}\\
\caption{Map $L$ defined as the composition of the maps described by Eqs. (\ref{Map}) and (\ref{R}).}\label{M5}
\end{center}
\end{figure}

Consider first the dynamical system $M: \mathcal{U} \to \mathcal{U}$, with phase space $\mathcal{U}:=[0,1]\times [0,1]$, defined by:
\be
\left(\begin{array}{c}
 x_{n+1} \\
 y_{n+1}
\end{array}\right)
 = M\left(\begin{array}{c}
 x_{n} \\
 y_{n}
 \end{array}\right)=
\left\{
 \begin{array}{l c}
   \left(\begin{array}{c}
           \dfrac{1}{2\ell}x_{n}+\dfrac{1}{2}\\ \\
           \dfrac{1}{2} y_{n} +\dfrac{1}{2}
         \end{array}\right) & \quad \text{for $0\leq x < \ell$}\\
 & \\
       \left(\begin{array}{c}
           \dfrac{x_{n}}{1-2\ell}-\dfrac{\ell}{1-2\ell} \\ \\
           (1-2\ell) y_{n} + 2\ell
         \end{array}\right) & \quad \text{for $\ell \leq x <
\frac{1}{2}$} \\
 & \\
       \left(\begin{array}{c}
           2x_{n}-\dfrac{1}{2} \\ \\
           \dfrac{1}{2} y_{n}
         \end{array}\right) & \quad \text{for $\frac{1}{2} \leq x < \frac{3}{4}$} \\
 & \\
       \left(\begin{array}{c}
           2x_{n}-\dfrac{3}{2} \\ \\
           2\ell y_{n}
         \end{array}\right) & \quad \text{for $\frac{3}{4} \leq x \leq 1$}
 \end{array}\right. \quad , \label{Map}
\ee
which is time reversal invariant in the sense that there exists an involution $G$ such that
\be
MGM = G \label{rev1} \quad .
\ee
Let $J_M (\vec{x})$ be the jacobian determinant of $M$ at $\vec{x}$, and
\be
\Lambda^{(M)} (\vec{x}) = \ln J_{M}^{-1}(\vec{x}) \label{invol}
\ee
be the phase space contraction rate. Then, $M$ corresponds to the  ``equilibrium'' version of a more general model introduced in 
Refs. \cite{CR,CKDR}, in the sense that its ensemble average phase space contraction rate vanishes.
Averaging in time over a trajectory of $n$-steps, yields:
\be
\overline{\Lambda}_n^{(M)} (\vec{x})  = \frac{1}{n}\sum_{k=0}^{n-1} \Lambda^{(M)}(M^k \vec{x}) \quad .
\ee
Next, let us introduce the rotation defined by the map $R$:
\be
\left(\begin{array}{c}
 x' \\
 y'
\end{array}\right)
 = R\left(\begin{array}{c}
 x \\
 y
\end{array}\right)= 
\left(\begin{array}{c}
 1-y \\
 x
\end{array}\right)  \quad .\label{R}
\ee
The composition $L=R\circ M$, illustrated by Fig. \ref{M5}, is a map whose phase space contraction rate $ \Lambda $ equals that of $M$, i.e.: 
$\Lambda(\vec{x}) = \Lambda^{(M)} (\vec{x})$. This does not mean that the time averages $\overline{\Lambda}_n(\vec{x})$ 
and $\overline{\Lambda}_n ^{(M)}(\vec{x})$ are also equal, because the trajectories of $L$ differ from those of $M$.

\section{The dynamics}
\label{sec:dyn}
\begin{figure}[htbp!]
\begin{center}
\raisebox{-1mm}{\includegraphics[width=0.41\textwidth]{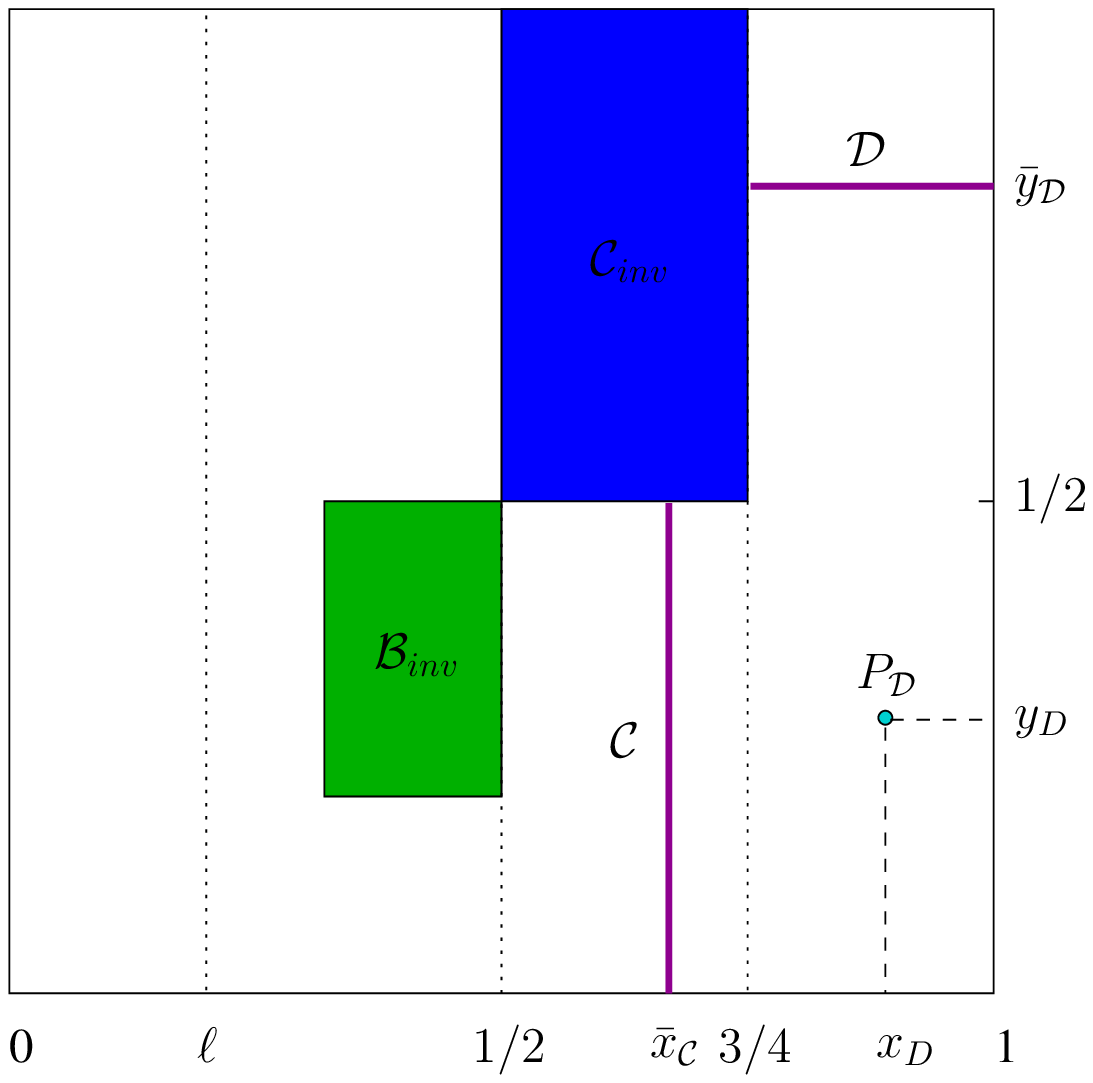}}
\includegraphics[width=0.41\textwidth]{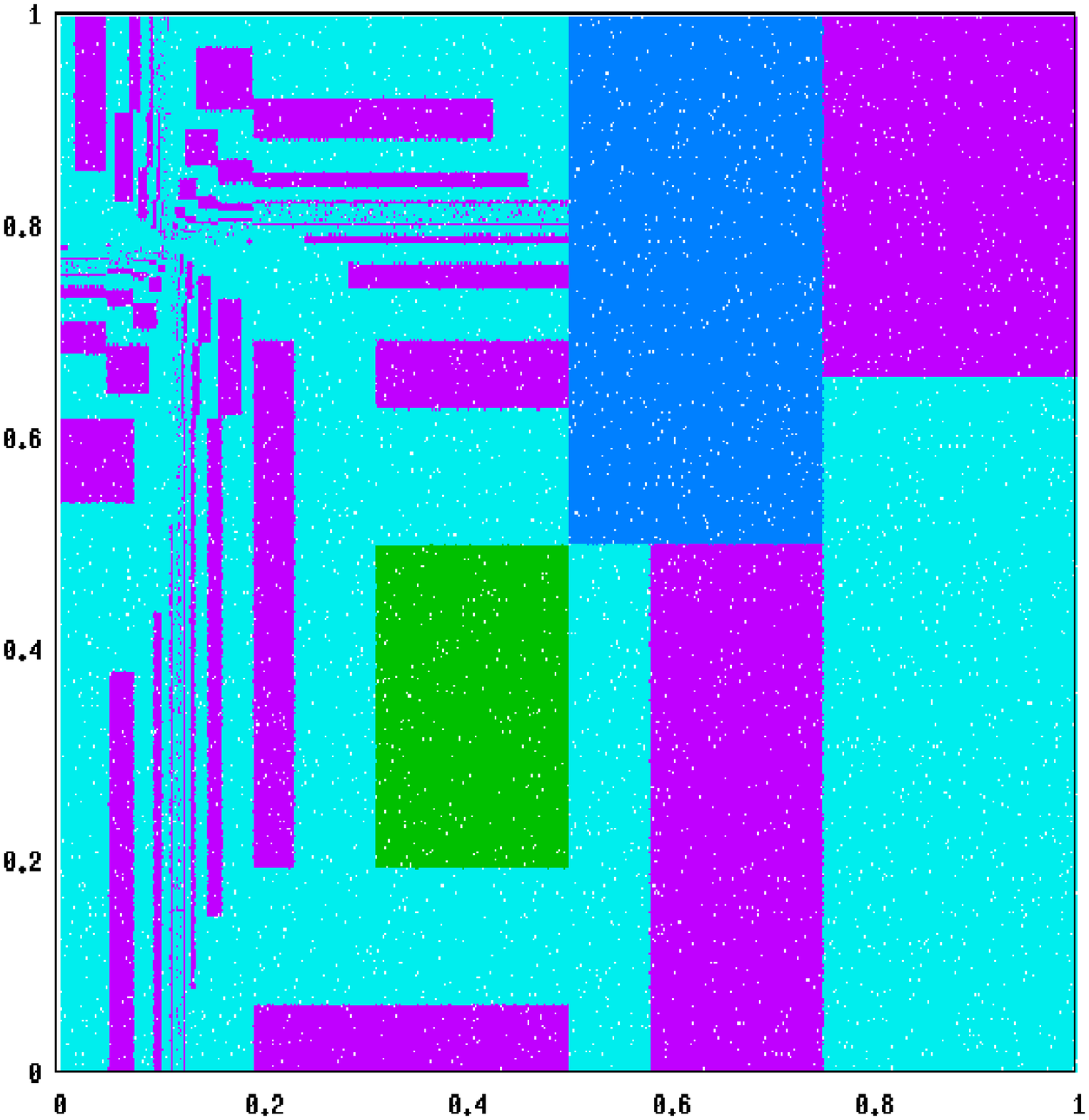}
\caption{\textit{Left panel}: the attractors of the map $L$. \textit{Right panel}: the corresponding basins of attraction (shown is the case with $\ell = 0.19$). Points in the area coloured in turquoise will converge in the steady state to the attractor $P_\mathcal{D}$, while the points lying in the purple regions will collapse to the two orthogonal lines referred to, in the left panel, as $\mathcal{C}$ and $\mathcal{D}$. Those orbits are only possible if the parameter $\ell$ satisfies $\frac{1}{8} \leq \ell \leq \frac{1}{4}$. Finally the green and blue regions are invariant regions which remain unchanged under the dynamics, and coincide with the regions $B_{inv}$ and $C_{inv}$ portrayed on the left panel.}
\label{fig:attr}
\end{center}
\end{figure}

Let us introduce $\Phi = -\ln(|J_D|)$, where $J_D$ is the Jacobian of the transformation $L$ in the region $D$ (i.e. $x \in [\frac{3}{4},1]$). 
One has $|J_D|=|J_A|^{-1}$, where $J_A$ is the the Jacobian of the transformation $L$ in the region $A$ (i.e. $x \in [0,\ell)$). 
Hence, we can write :

\bea
\Lambda (x,y) = \Lambda (x) = \left\{
\begin{array}{c}
- \Phi  \quad  \text{for $ 0 \leq x < \ell $} \nonumber\\
0 \quad  \text {for $\ell \leq x < \frac{1}{2}$} \nonumber \\
0 \quad  \text{for $ \frac{1}{2} \leq x < \frac{3}{4}$} \nonumber\\
\Phi \quad \text{for $ \frac{3}{4}\leq x \leq 1 $} \nonumber
\end{array} \right. \quad .
\eea
Figure \ref{fig:attr} (on the left) shows the invariant sets of the map $L$. They consist of:
\begin{itemize}
\item
two invariant regions (coloured in green and dark blue), called $\mathcal{B}_{inv}$ and $\mathcal{C}_{inv}$, which are characterized by zero Lyapunov exponents;
\item
a fixed point $P_\mathcal{D}$, characterized by two negative Lyapunov exponents;
\item
two orthogonal lines (coloured in purple), labelled as $\mathcal{C}$ and $\mathcal{D}$, along which we have one vanishing exponent. At every point of these lines, the orthogonal directions correspond to one negative Lyapunov exponent.
\end{itemize}

More fixed points and cycles are present, which have not been evidenced in Fig. \ref{fig:attr} since they are not attractors: a repulsive fixed point $P_\mathcal{A}$ in the region $A$, characterized by two positive Lyapunov exponents, and a repulsive 
(hyperbolic) cycle, referred to as $\mathcal{AB}$ periodic orbit, constituted by two points, one in the $A$ and the other in the $B$ region. 
The area coloured in turquoise, on the right panel of Fig. \ref{fig:attr}, converges to the attracting fixed point $P_\mathcal{D}$, while the points lying in the purple regions collapse on cycles constituted by $4$ points lying on the two orthogonal lines labelled by $\mathcal{C}$ and $\mathcal{D}$. We refer to these cycles as $\mathcal{C}\mathcal{D}\mathcal{C}\mathcal{D}$ cycles. 
The central invariant regions, coloured in green and blue, are instead the same as the invariant coloured regions $B_{inv}$ and $C_{inv}$ on the left. Indeed, points lying in this two areas will move on a period-4 cycle confined within the same borders.
These regions are invariant, as it is easy to show analytically: applying four times the evolution operator to any point of such regions, the evolution remains in the same area and eventually returns to the starting point.
While the height and width of the image of the invariant region $B_{inv}$, under the mapping $L$, depend on $\ell$, the height and width of the image of the invariant region $C_{inv}$ do not. \\
Let us compute the corresponding Lyapunov exponents. To this aim, we denote by $DL(x,y)$ the Jacobian matrix of the map $L$ and notice that, in the region $A$, that matrix takes the form:
\be
DL\bigg|_A = \begin{bmatrix} 0 & -\frac{1}{2} \\ \frac{1}{2 \ell} & 0 \end{bmatrix} \quad .
\ee
Therefore, the Lyapunov exponents for the fixed point $P_\mathcal{A}$, with coordinates
$$
P_\mathcal{A} = (x_A,y_A) = \left(\dfrac{\ell}{1+4\ell},\dfrac{1+2 \ell}{1+4\ell}\right) \quad ,
$$
along the horizontal and  vertical directions, $\zl_x$ and $\zl_y$ respectively, are given by
\bea
\lambda_x  (x_A, y_A) &=& 
%\lim_{n \to \infty} \frac{1}{n} \ln \left| \left(DL\bigg|_A \right)^n \right| =  
\lim_{n \to \infty} \frac{1}{n} \ln \left[ \left( \frac{1}{2 \ell}\right)^\frac{n}{2} \left(\frac{1}{2}\right)^\frac{n}{2} \right] = \nonumber \\
&=& \frac{1}{2} \ln \left( \frac{1}{4 \ell}\right) = \lambda_y (x_A,y_A) \quad ,
\eea
which is positive $ \forall \ell \in (0, \frac{1}{4})$. 
It follows that the fixed point $P_\mathcal{A}$ is a repeller and that two phase space points, initially in its neighbourhood, will diverge exponentially fast under the dynamics. 

In the contracting region $D$, the fixed point is given by
$$
P_\mathcal{D} = (x_D,y_D) = \left(\dfrac{1+3\ell}{1+4 \ell},\dfrac{1}{2(1+4 \ell)} \right) \quad ,
$$ 
the corresponding jacobian matrix can be written as
\be
DL\bigg|_D = \begin{bmatrix} 0 & -2 \ell \\ 2 & 0 \end{bmatrix} \quad ,
\ee
and the Lyapunov exponents thus read
\bea
\lambda_x  (x_D, y_D) =\frac{1}{2} \ln  \left( 4 \ell \right) = \lambda_y (x_D,y_D) \quad .
\eea
Note that $\lambda_x  (x_D, y_D)$ and $\lambda_y (x_D,y_D)$ are negative $\forall \ell \in (0, \frac{1}{4})$, and opposite to $\lambda_x (x_A,y_A)$ and $\lambda_y (x_A,y_A)$.

For $ \ell \geq \frac{1}{8}$ it is possible to show analytically that a pair of new conjugate trajectories, corresponding to the periodic orbits $\mathcal{AB}$ and $\mathcal{CDCD}$, arise. 
One of them jumps from the neutral area $B$ to the expanding area $A$, while the latter is alternatively stepping from the neutral region $C$ to the contracting region $D$. Hence, the corresponding value of $\overline{\Lambda}_n$ attains, for large $n$, respectively, the asymptotic values $-\Phi/2$ and $\Phi/2$.
\begin{figure}[htbp!]
\begin{center}	
\includegraphics[scale=0.25]{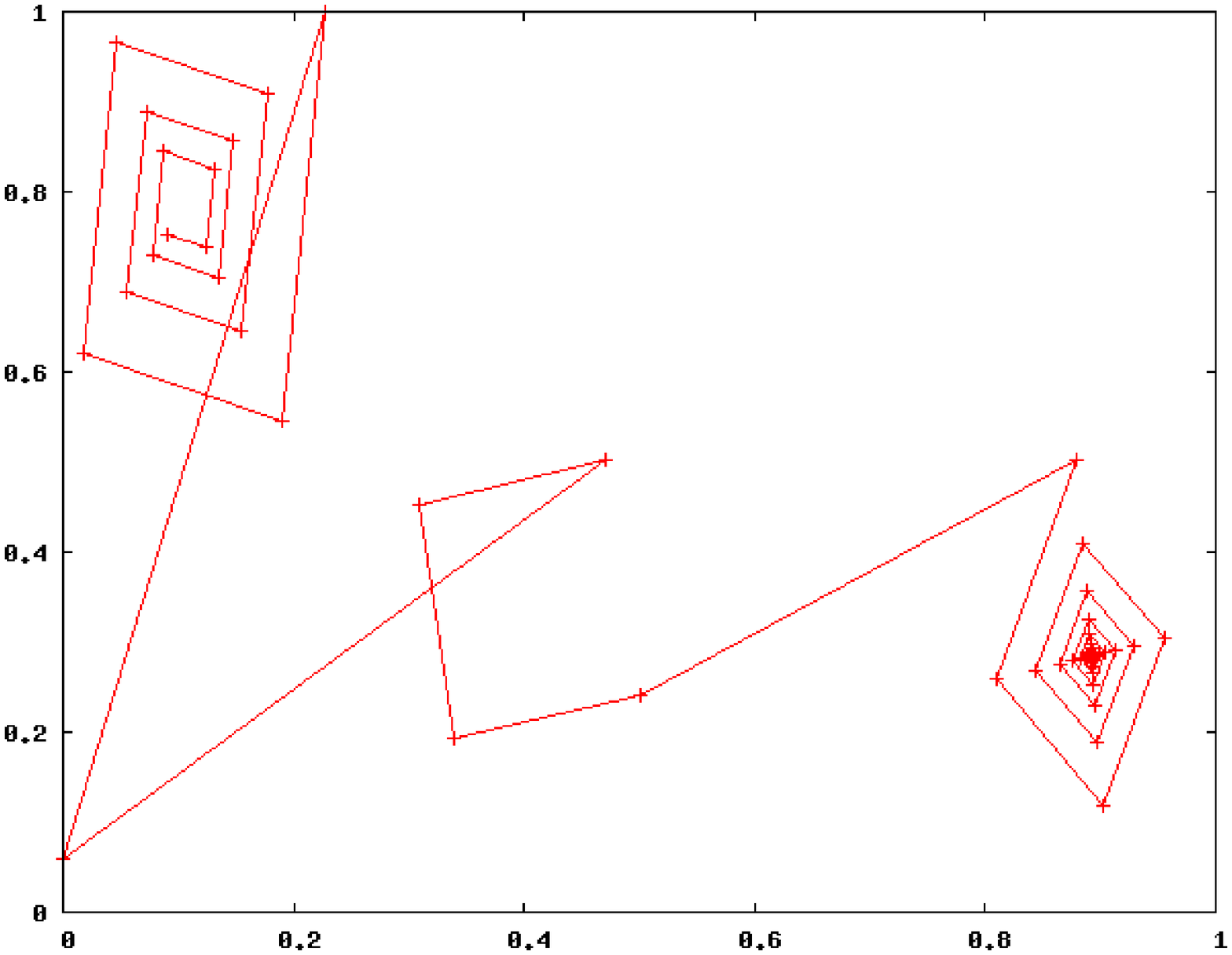}
\includegraphics[scale=0.25]{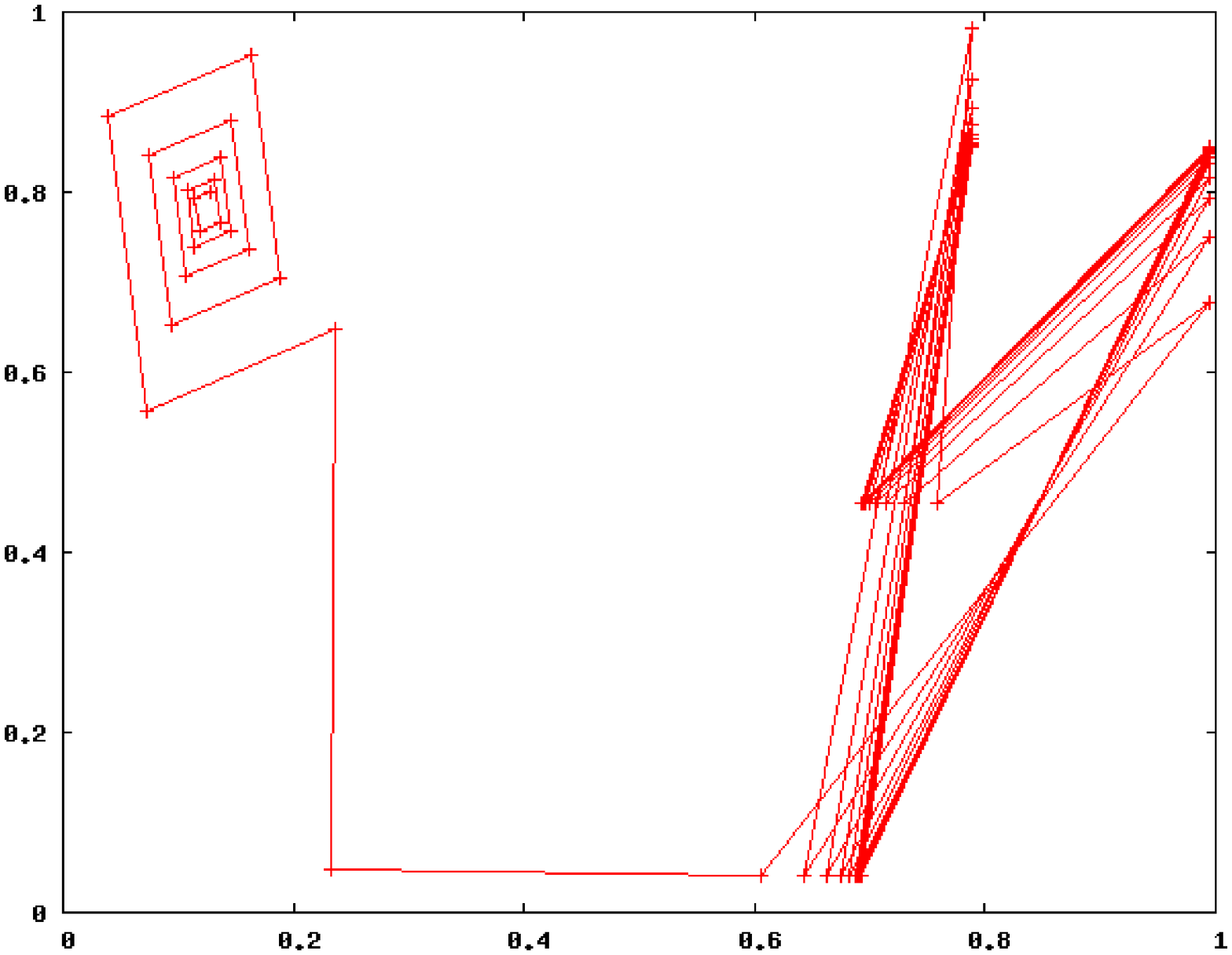}
\caption{ {\itshape Left Panel:} a sample trajectory starting in the neighbourhood of $P_\mathcal{A}$ and evolving towards $P_\mathcal{D}$. {\itshape Right panel:} a sample trajectory evolving towards the $\mathcal{CDCD}$ orbit.}
\label{traj}
\end{center}
\end{figure}

The periodic orbit $\mathcal{CDCD}$ lies on the two orthogonal lines, labeled by $\mathcal{C}$ and $\mathcal{D}$ in Fig. \ref{fig:attr}, and identified, respectively, by the coordinates 
$$
\bar{x}_{\mathcal{C}}=\dfrac{1+\ell}{1+4\ell} \quad \text{and} \quad \bar{y}_{\mathcal{D}}=\dfrac{3}{2(1+4\ell)} \quad .
$$
Every point lying on those lines is a point of a period-4 cycle which jumps from the vertical to the horizontal line, and viceversa.
The corresponding Lyapunov exponents vanish in the direction of the lines $\mathcal{C}$ and $\mathcal{D}$, and are negative along the orthogonal direction. 
For instance, taking a point on the line $\mathcal{C}$, one obtains:
\be
\lambda_x  (\bar{x}_{\mathcal{C}}, y) = \lim_{n \to \infty} \frac{1}{n} \ln \left[ \left( 2 \ell \right)^\frac{n}{2} \left(2\right)^\frac{n}{2} \right] = \frac{1}{2} \ln  \left( 4 \ell \right) \quad ,
\ee
which is negative $\forall \ell \in (0,\frac{1}{4})$, and
\be
\lambda_y  (\bar{x}_{\mathcal{C}}, y) = \lim_{n \to \infty} \frac{1}{n} \ln \left[ \left( 2 \right)^\frac{n}{2} \left(\frac{1}{2}\right)^\frac{n}{2} \right] = 0 \quad .
\ee
Analogously, we have $\lambda_x(x,\bar{y}_{\mathcal{D}}) = 0$ and $\lambda_y(x,\bar{y}_{\mathcal{D}}) < 0$, where by $(\bar{x}_{\mathcal{C}}, y)$ we mean any point lying on the line $\mathcal{C}$ and by $(x,\bar{y}_{\mathcal{D}})$ any point belonging to the line $\mathcal{D}$.

The conjugate trajectory with respect to the orbit $\mathcal{CDCD}$ is the period-2 orbit $\mathcal{AB}$, which is constituted by two points, one in the region $A$, with coordinates
$$
(\bar{x}_{\mathcal{A}},\bar{y}_{\mathcal{A}})= \left(\ell (1-2 \ell),\dfrac{1-2\ell}{3-4 \ell} \right) \quad ,
$$
and the other in the region $B$, with coordinates
$$
(\bar{x}_{\mathcal{B}},\bar{y}_{\mathcal{B}})= \left(\dfrac{1-\ell}{3-4 \ell},1-\ell \right) \quad .
$$
The evaluation of the Lyapunov exponents for the point belonging to $\mathcal{AB}$ and lying in the region $A$ (but the same result also holds for the point located in $B$) yields
\bea
\lambda_x  (\bar{x}_{\mathcal{A}}, \bar{y}_{\mathcal{A}}) &=& \lim_{n \to \infty} \frac{1}{n} \ln \left[ \left(1 - 2 \ell \right)^\frac{n}{2} \left(\frac{1}{2 \ell}\right)^\frac{n}{2} \right] = \nonumber \\
&=& \frac{1}{2} \ln  \left( \frac{1-2\ell}{2\ell} \right) \quad ,
\eea
which is positive for each $\ell \in (0, \frac{1}{4})$, and
\bea
\lambda_y  (\bar{x}_{\mathcal{A}}, \bar{y}_{\mathcal{A}}) &=& \lim_{n \to \infty} \frac{1}{n} \ln \left[ \left(\frac{1}{2}\right)^\frac{n}{2} \left(\frac{1}{1- 2 \ell}\right)^\frac{n}{2} \right] = \nonumber \\
&=& \frac{1}{2} \ln  \left( \frac{1}{2(1-2\ell)} \right) \quad ,
\eea
which is negative for each $\ell \in [0, \frac{1}{4})$. 
From the computation of the Lyapunov exponents, it emerges that the points on the orbits $\mathcal{AB}$ are hyperbolic repellers, thus the dynamics escapes from such cycles. 

Moreover, we have $(\lambda_x + \lambda_y)_{\mathcal{AB}} = - (\lambda_x + \lambda_y)_{\mathcal{CDCD}}$, as expected for conjugate trajectories, which must produce opposite values of $\bar{\Lambda}_n$, but differently from the standard time reversal invariant cases, the orbit conjugated to a $\mathcal{CDCD}$ orbit has period 2 rather than 4. 
For the validity of certain symmetries such as the FR, this does not necessarily constitute a difficulty, since all trajectories come in pairs with opposite average phase space contraction.

\begin{figure}[htbp!]
\begin{center}	
\includegraphics[scale=0.25]{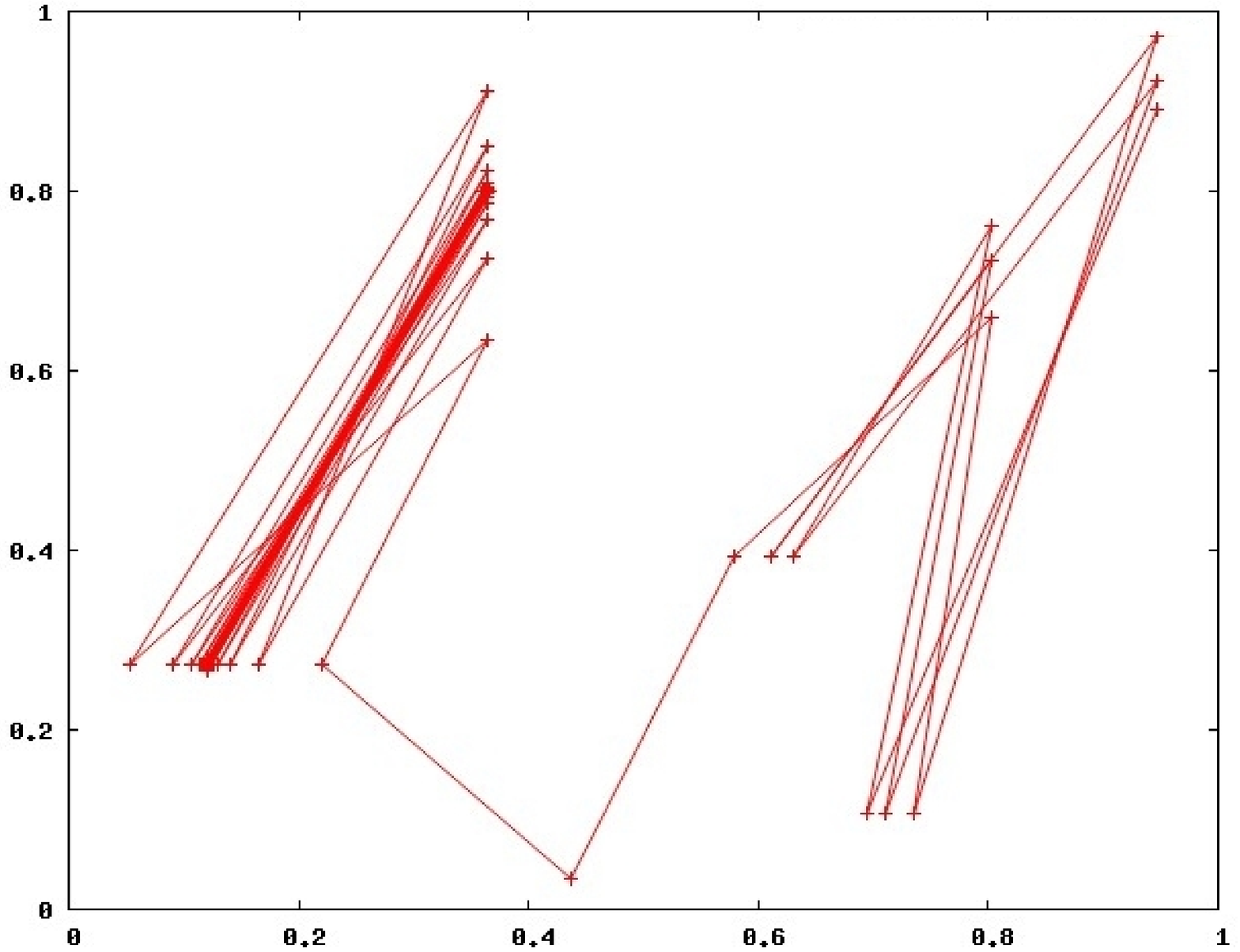}
\includegraphics[scale=0.25]{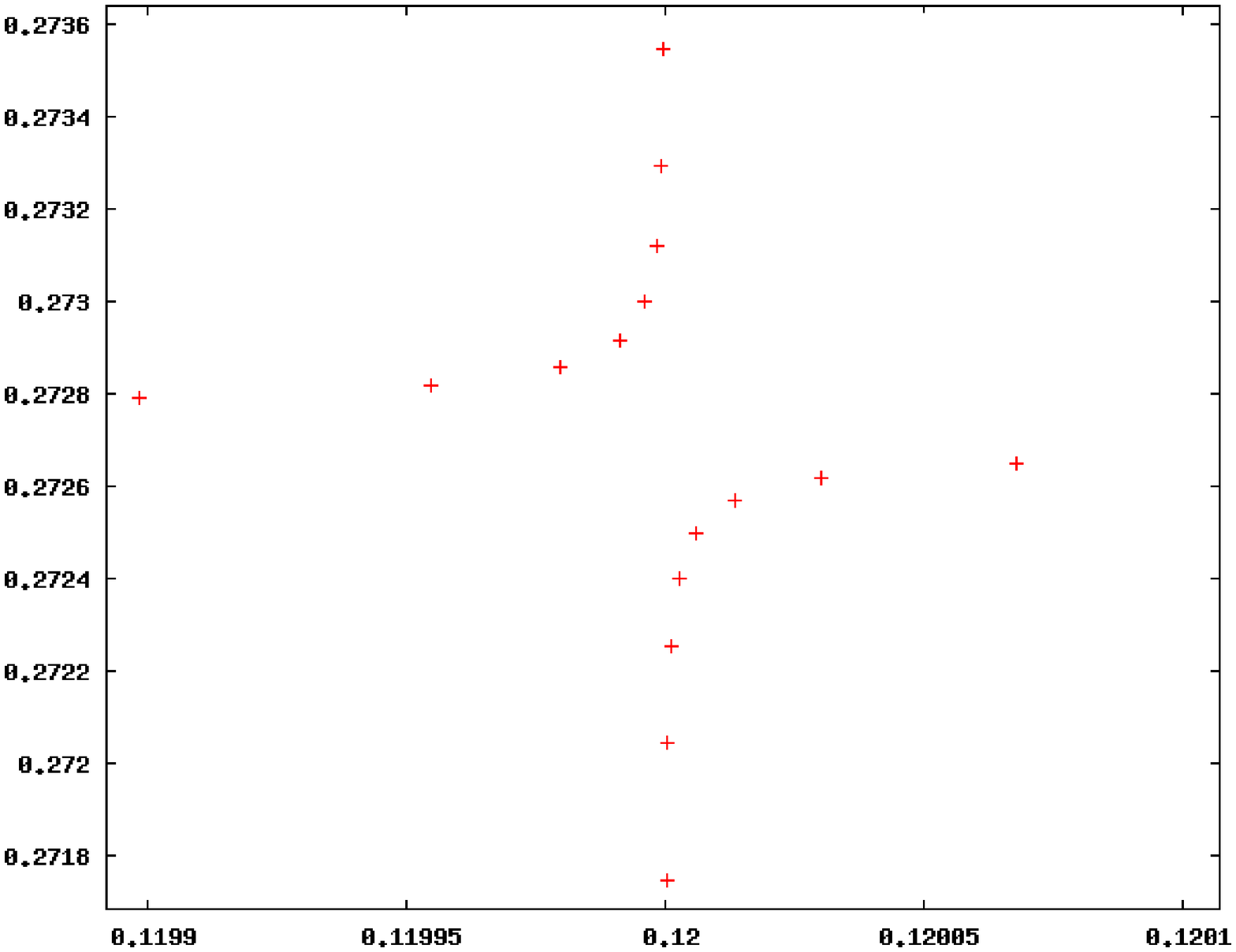}
\caption{{\itshape Left panel:} a sample of a trajectory escaping from the$\mathcal{AB}$ repulsive cycle and collapsing on the attractive $\mathcal{CDCD}$ attractive cycle. In the stationary state, the $\mathcal{AB}$ cycle and the $\mathcal{CDCD}$ cycle are conjugate trajectories, i.e. the sum of the Lyapunov exponents in the $x$ and $y$ components are opposite, $(\lambda_x + \lambda_y )_{\mathcal{AB}} = - (\lambda_x + \lambda_y )_{\mathcal{CDCD}}$. {\itshape Right panel:} detail of a trajectory in the surrounding of the hyperbolic point $P_\mathcal{A}$.}\label{ant}
\end{center}
\end{figure}

In Figure \ref{traj}, two trajectories starting in the surrounding of the repulsive fixed point $P_\mathcal{A}$ and attracted by different attractors (the fixed point $P_\mathcal{D}$ on the left and the $\mathcal{CDCD}$ orbit on the right) are shown. This reveals a fractal structure for the basins of attraction, see Fig. \ref{basins}.
However in the surrounding of the hyperbolic points there are no chaotic trajectories.
In Figure \ref{ant}, we illustrate one of these non-chaotic trajectories which starts around the unstable point in region B belonging to the periodic orbit $\mathcal{AB}$.

Although it has not been possible to find an involution satisfying the standard equation for the time reversal invariant dynamics, it has been numerically tested that for any finite number of steps, it exists a nonvanishing probability of finding pairs of trajectories  producing opposite values for the phase space contraction time-average $\bar{\Lambda}_n$.
For example, starting from a uniform distribution, for any given accessible value $\bar{\Lambda}$ produced by iterating the map for a number $n$ of steps, the probabilities of finding the corresponding trajectories can be computed exactly as the sum of the areas 
of the points $(x,y)$ satisfying the relation: 
$$
\bar{\Lambda}_n(L^n(x,y))=\bar{\Lambda} \quad .
$$ 
In particular, we found that the trajectories yielding negative values of $\bar{\Lambda}_n$ originate 
from points located close enough to the fixed unstable point in the expanding region: the closer the initial point is 
to the repulsor, the longer the trajectory will explore the expanding region before eventually collapsing on the attractor 
(the more negative, hence, the resulting value of $\bar{\Lambda}_n$ will be). We remark that this description applies to 
the transient regime (i.e. finite length trajectories), which, at variance with the steady state regime, allows fluctuations.

\section{Basins of attraction}

\begin{figure}
\centering
\includegraphics[width=11.5cm]{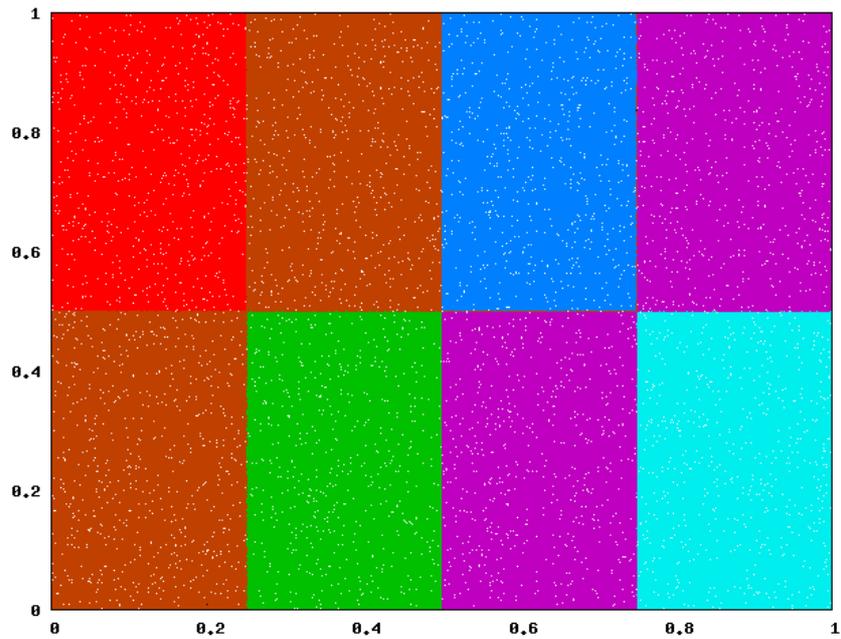}
\caption{The uniform ``equilibrium'' distribution corresponding to $\ell = 0.25$. No attractors are present, all the Lyapunov exponents vanish and each point belongs to a period-4 cycle, which remains within a same colour area. The phase space is fragmented in six invariant sets.}\label{eq}
\end{figure}

\begin{figure}[htpb!]
\centering
\includegraphics[width=11.5cm]{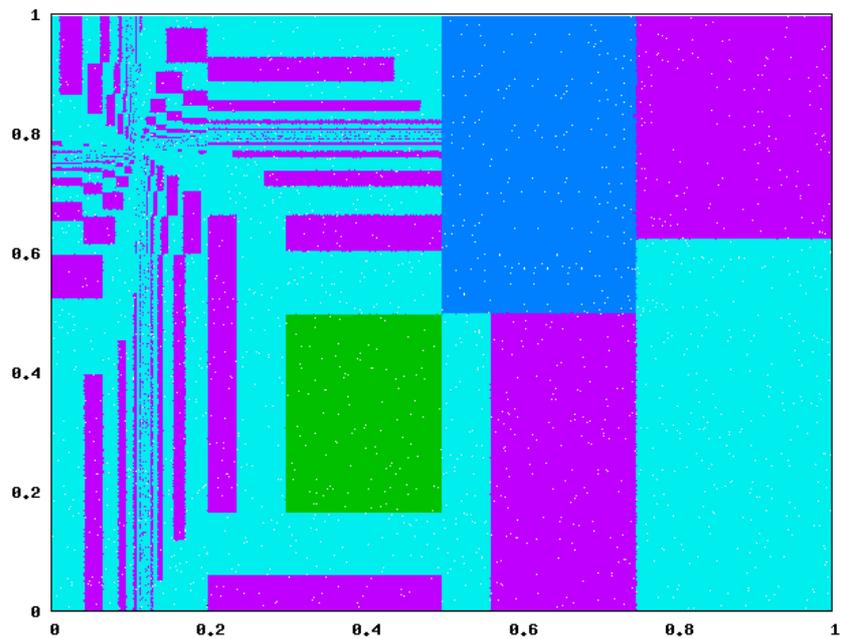}
\caption{Basins of attraction and invariant regions for the map illustrated in Fig. \ref{M5}, with $\ell=0.2$.}
\label{basins}
\end{figure}

We analyze the basins of attraction of the different attractors as functions of the parameter $\ell$.
In the equilibrium condition ($\ell=0.25$), shown in Fig. \ref{eq}, there are no attractors, all the Lyapunov exponents vanish and each point belongs to a period-$4$ cycle remaining within one of the six invariant regions.\\
Figure \ref{basins} shows the basins of attraction for $\ell = 0.2$: the green and blue sets correspond to the set of points which fall in the invariant regions $\mathcal{B}_{inv}$ and $\mathcal{C}_{inv}$ described in the previous section. The phase space area coloured in turquoise corresponds to the set of points whose dynamics collapse on the fixed point $P_\mathcal{D}$, whose coordinates, as shown in Sec. \ref{sec:dyn}, depend parametrically on $\ell$.
The remaining regions, coloured in purple, correspond to the basins of attraction of the $\mathcal{CDCD}$ orbits. These basins of attraction form complementary fractal sets stemming out of the unstable fixed point $P_\mathcal{A}$.

\begin{figure}[htbp!]
\center
\includegraphics[scale=0.21]{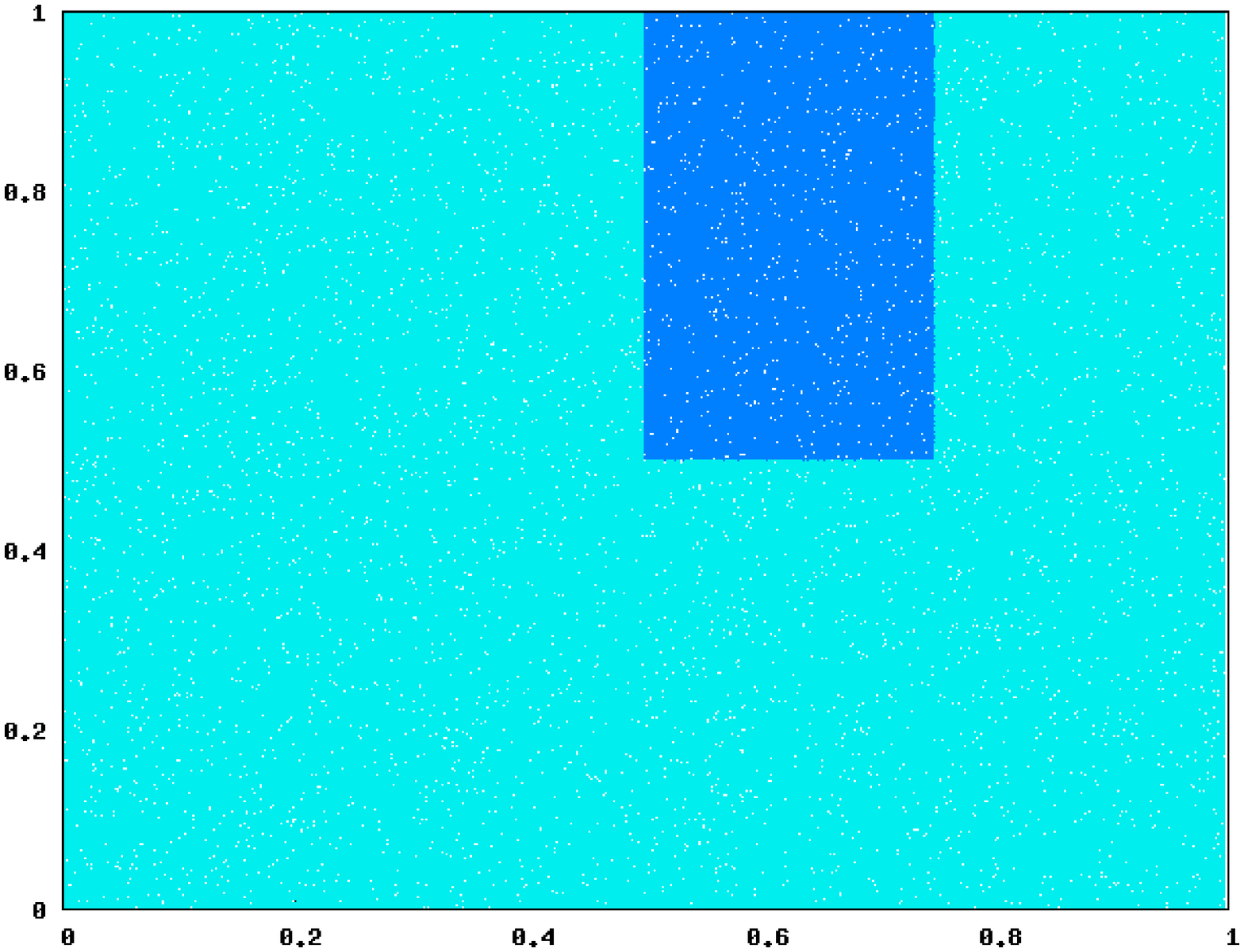}
\includegraphics[scale=0.21]{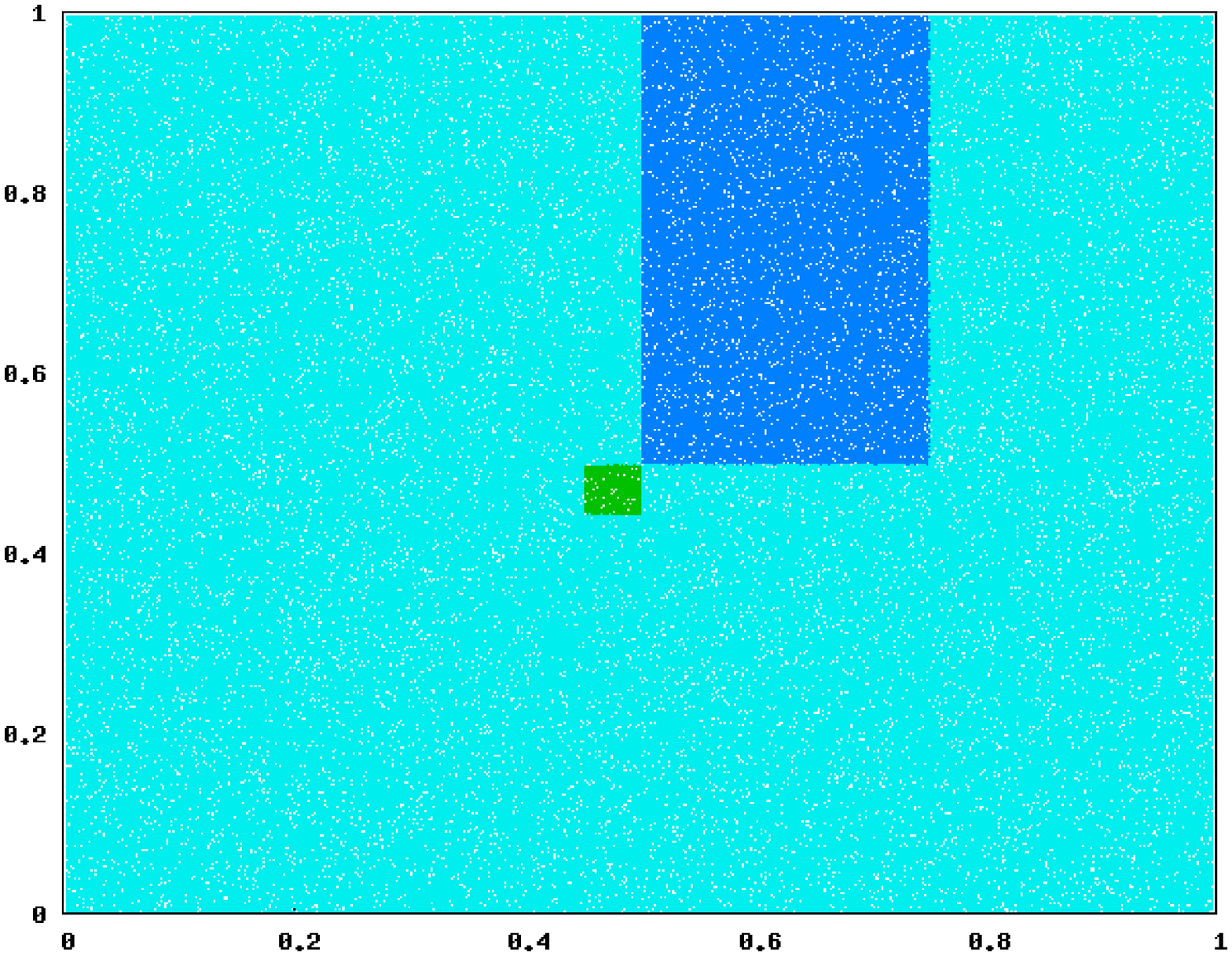}
\includegraphics[scale=0.21]{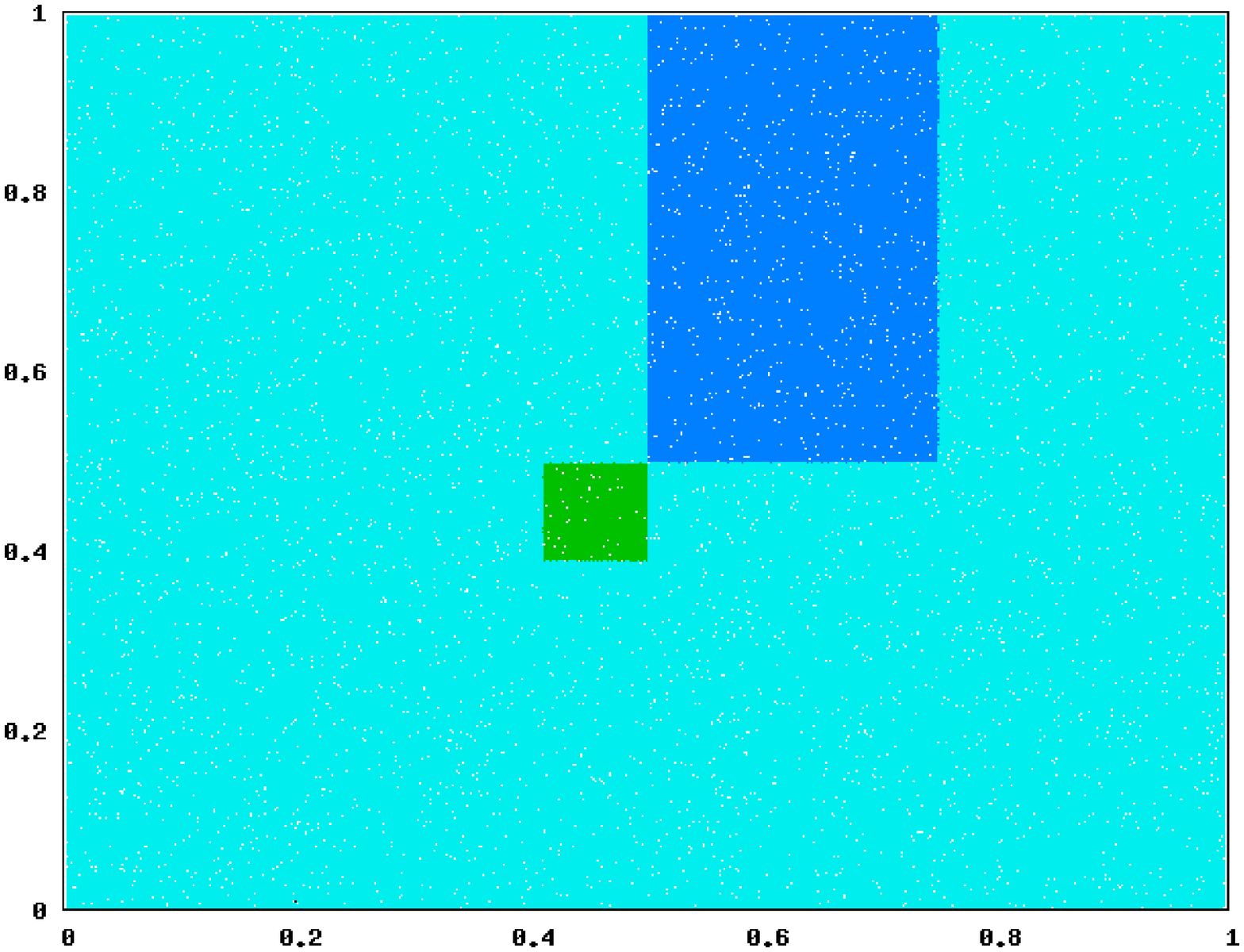}

\includegraphics[scale=0.21]{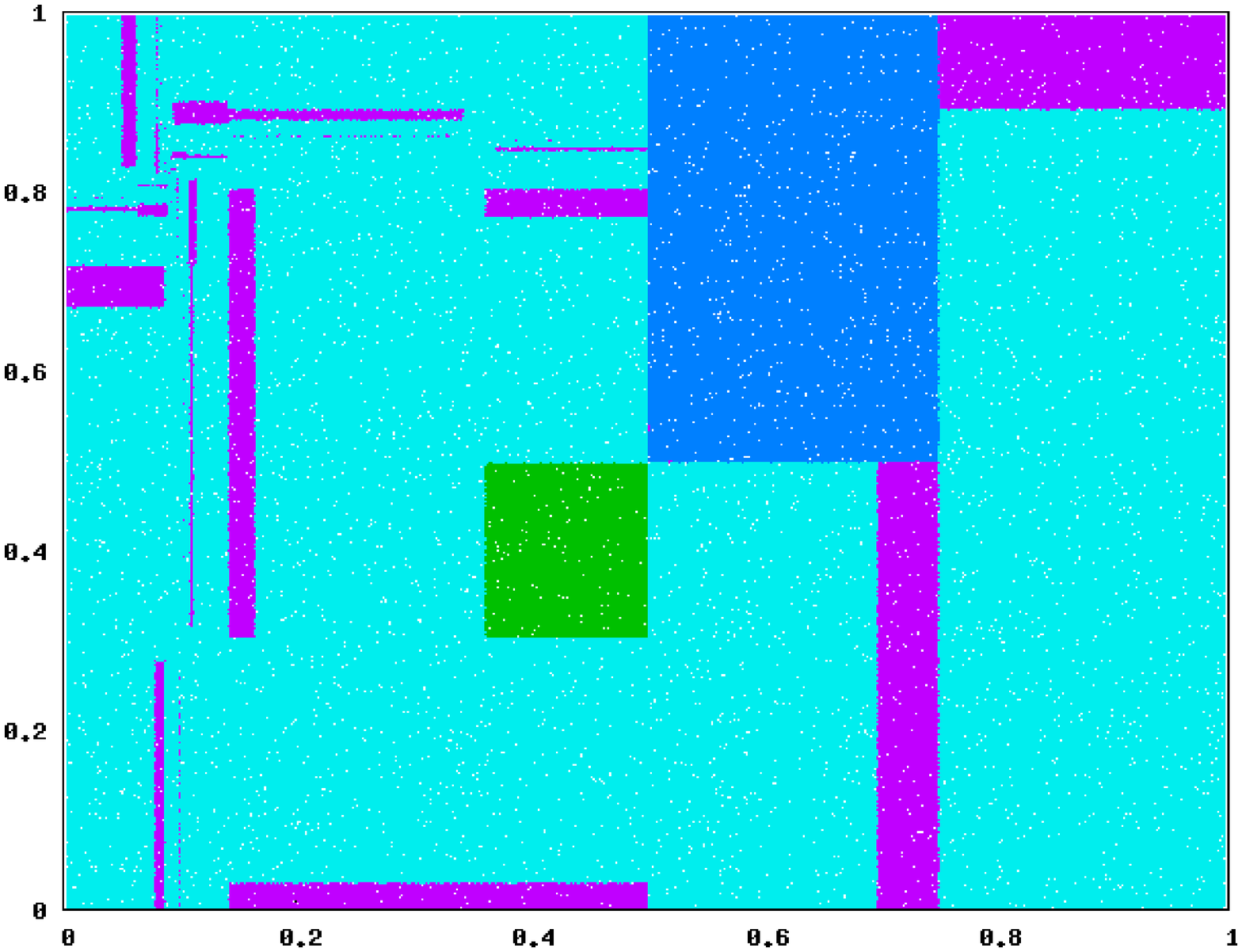}
\includegraphics[scale=0.21]{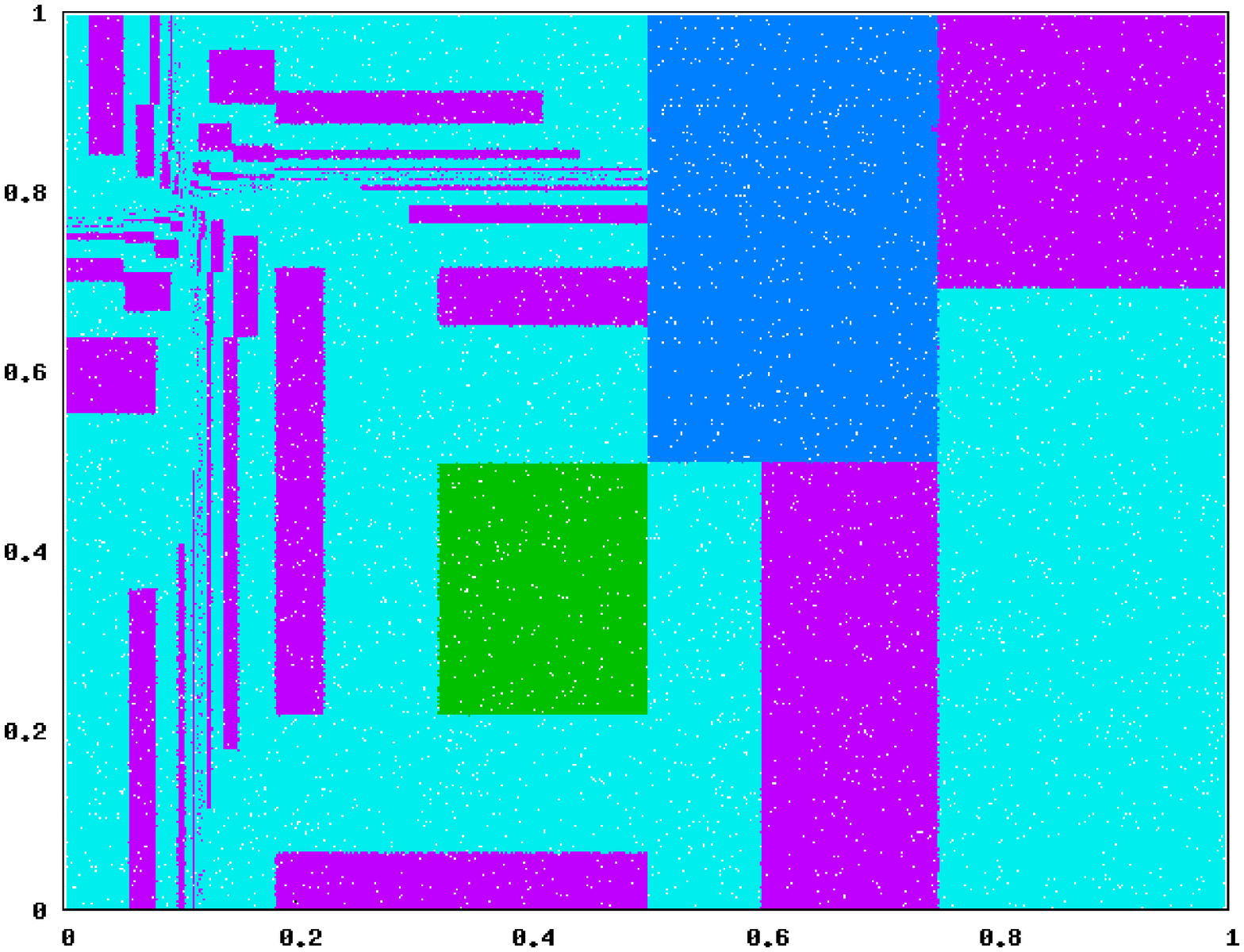}
\includegraphics[scale=0.21]{20.eps}

\includegraphics[scale=0.21]{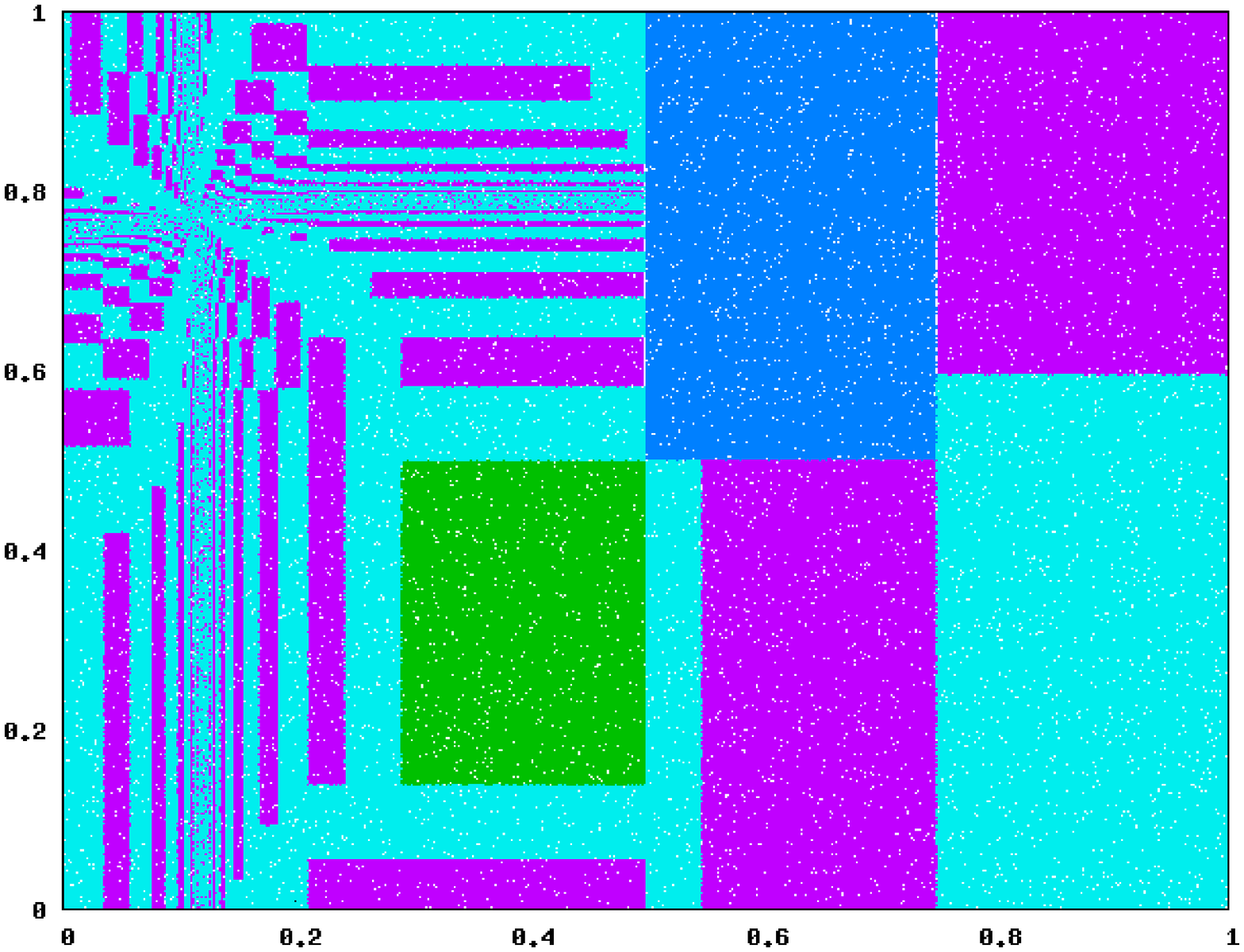}
\includegraphics[scale=0.21]{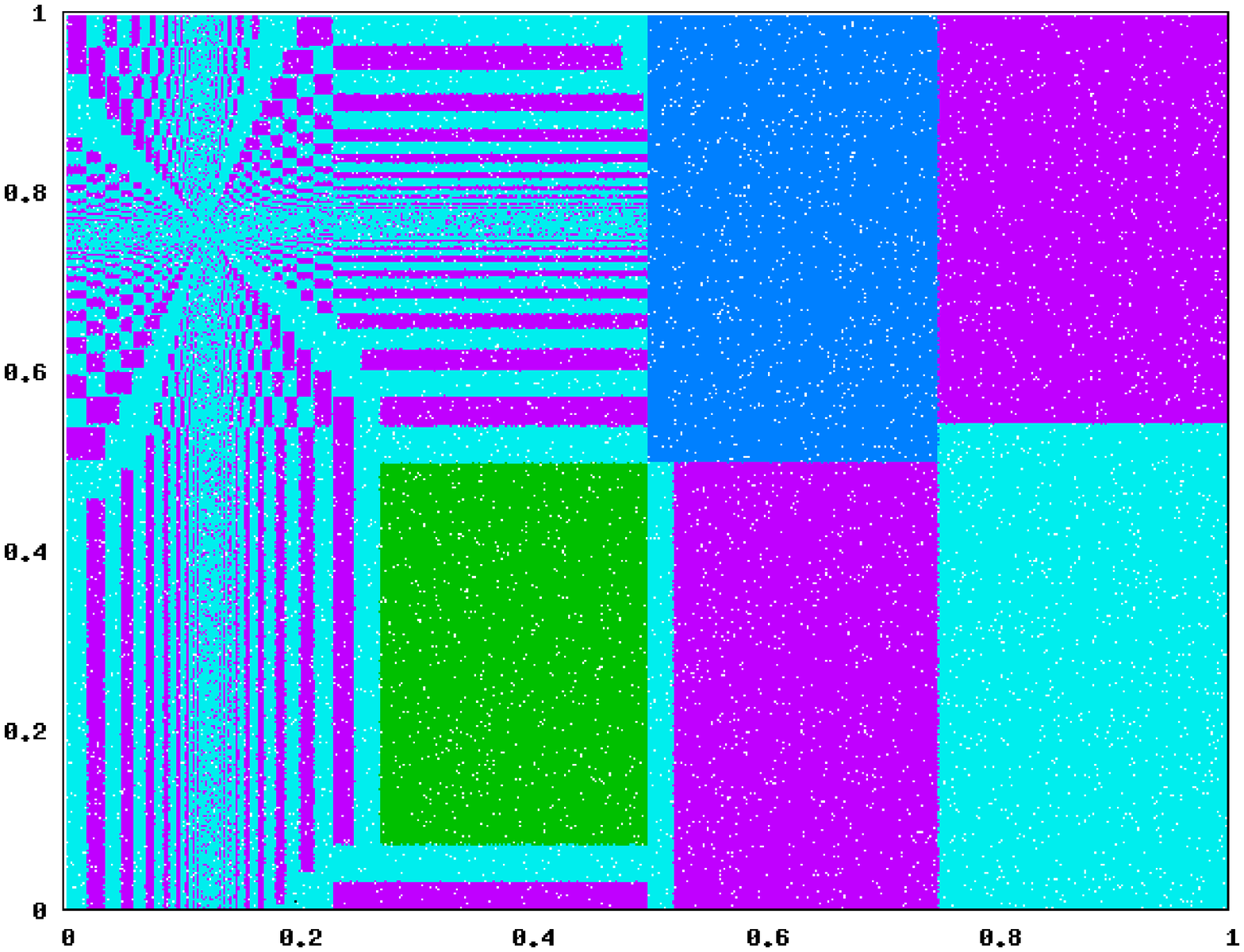}
\includegraphics[scale=0.21]{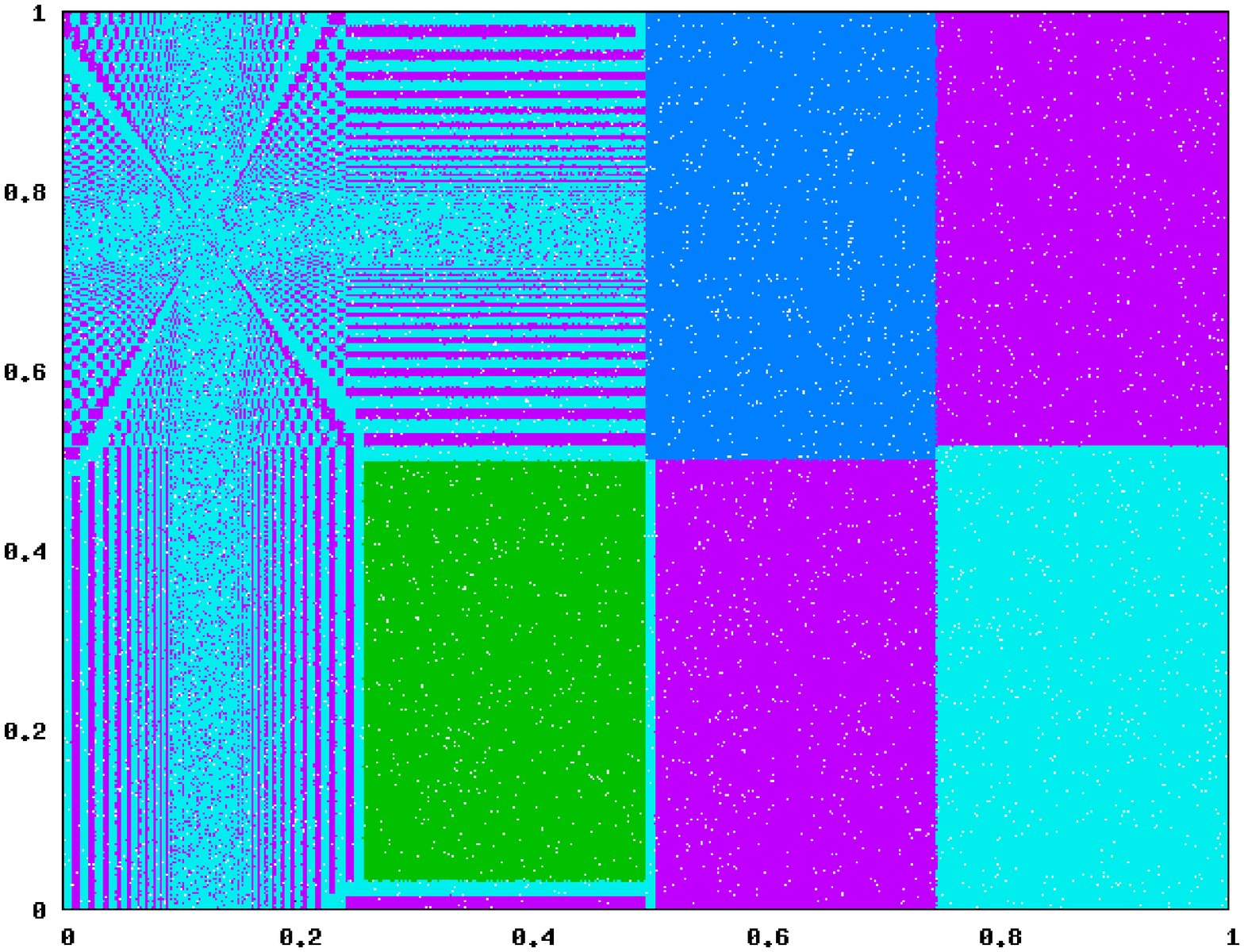}

\includegraphics[scale=0.21]{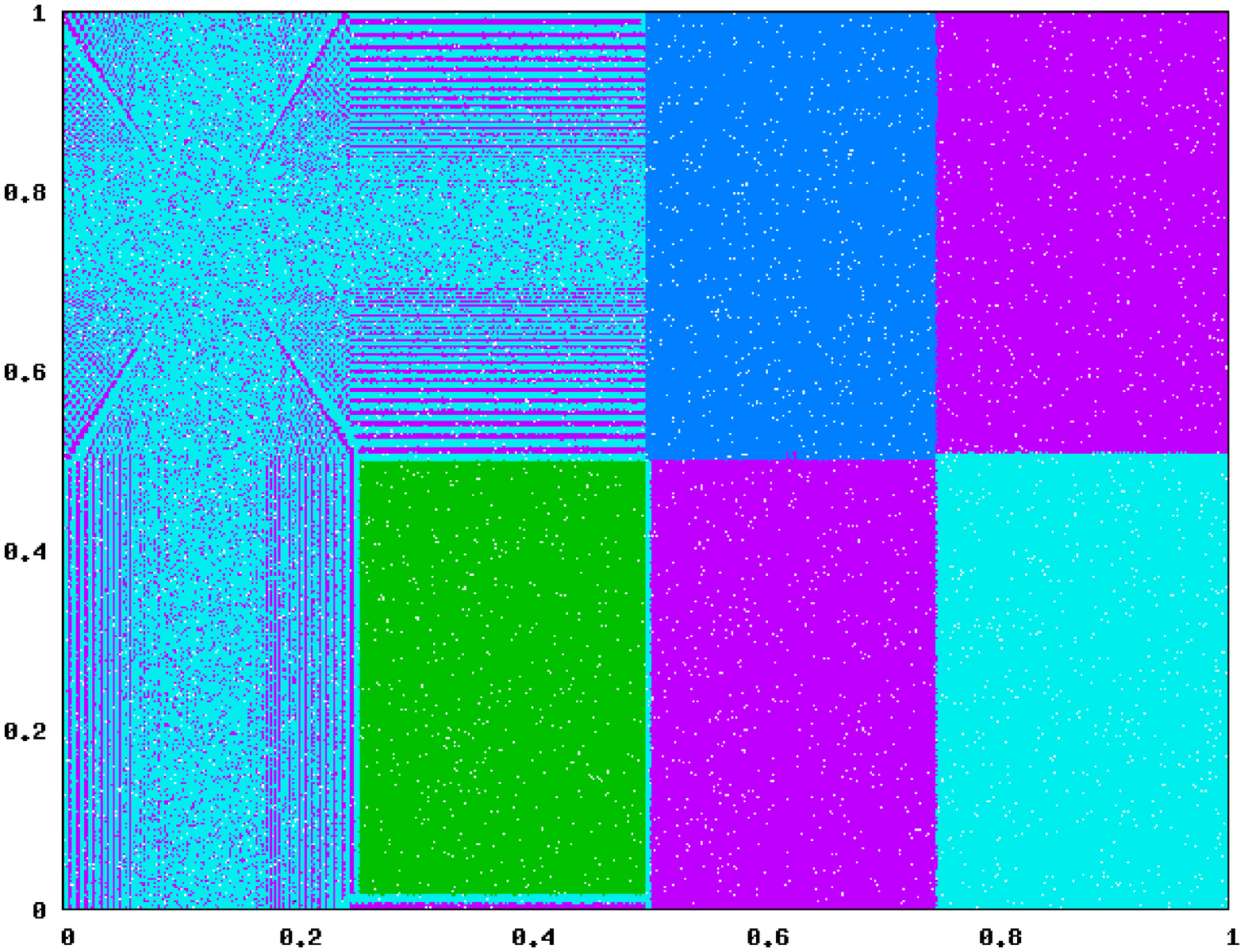}
\includegraphics[scale=0.21]{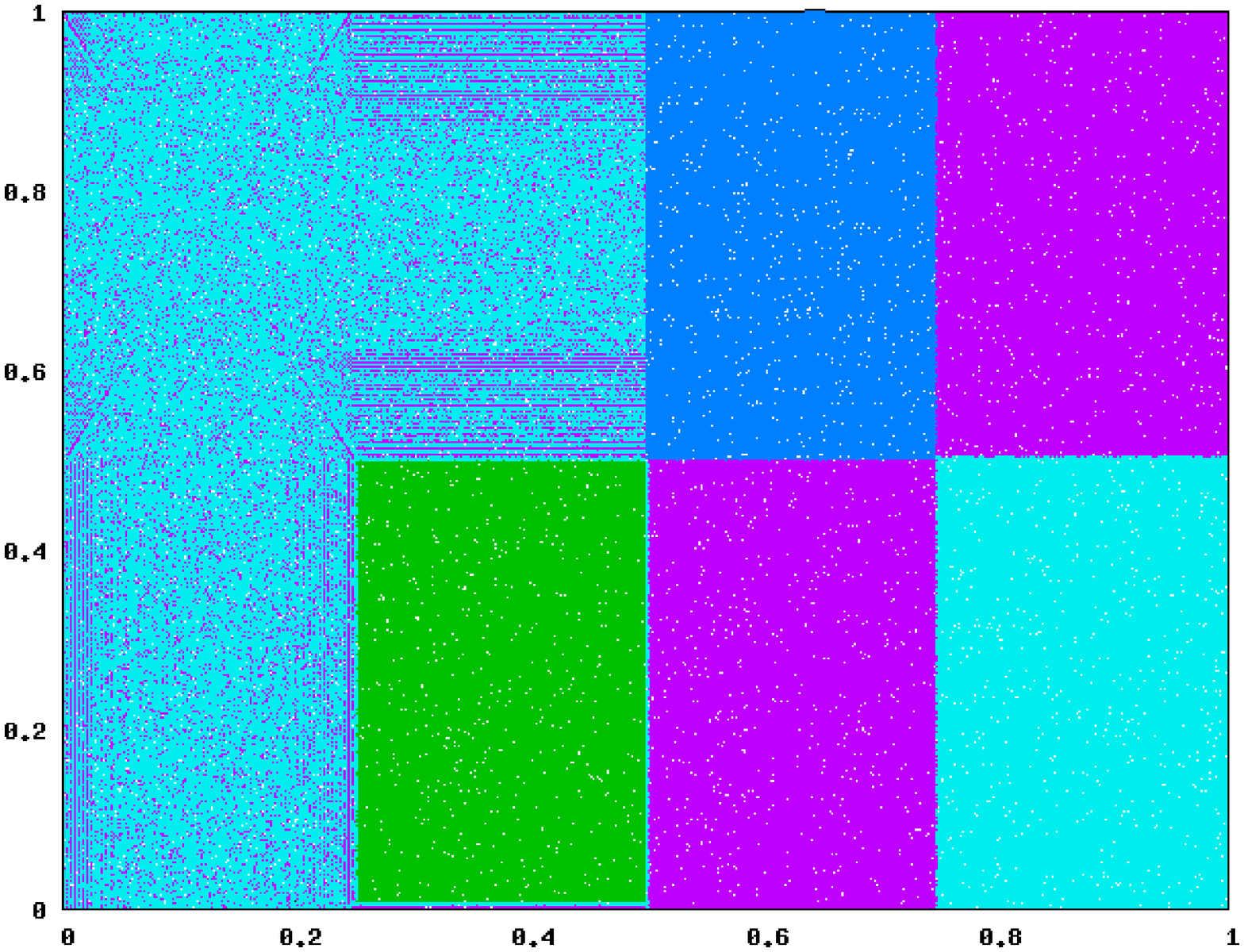}
\includegraphics[scale=0.21]{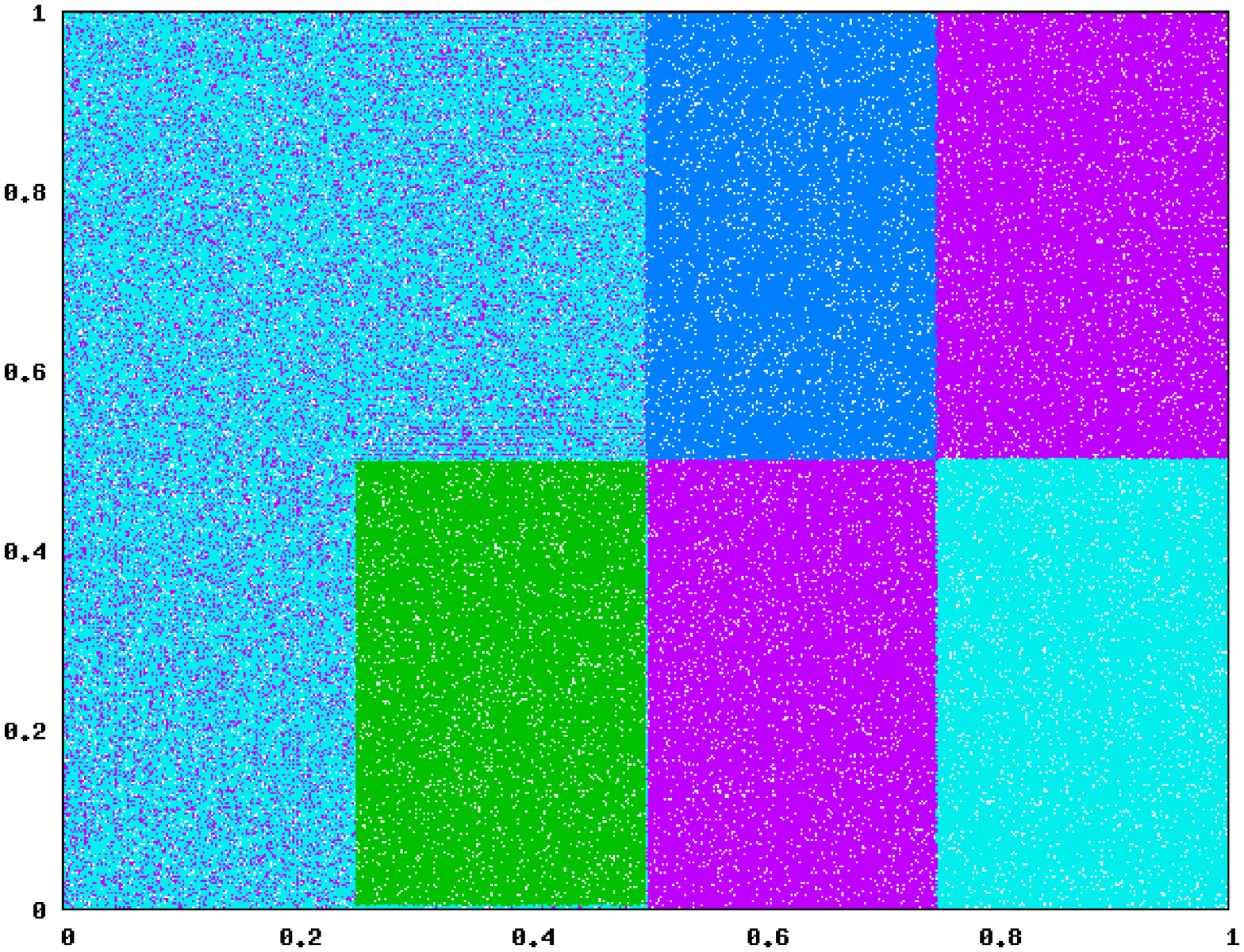}

\caption{Basins of attraction for increasing values of $\ell$. The blue and green regions correspond to invariant sets, the turquoise region is the basin of attraction of the fixed attractive point $P_\mathcal{D}$. From the second row, i.e. for $\ell\geq \frac{1}{8}$, the basins of attraction of $\mathcal{CDCD}$ orbits show up in purple, building a fractal around the unstable fixed point $P_\mathcal{A}$.
{\itshape From left to right, top to bottom:} $\ell=0.0001$, $\ell=0.05$, $\ell=0.09$; $\ell=0.14$, $\ell=0.18$ $\ell=0.20$; $\ell=0.021$ $\ell=0.23$, $\ell=0.242$; $\ell=0.246$ $\ell=0.248$, $\ell=0.249$.}
\label{bas}
\end{figure}

The set corresponding to the neutral invariant regions $\mathcal{B}_{inv}$ depends parametrically on $\ell$ and is expressed by:
$$ 
\frac{1}{2(1-2 \ell)} \leq x \leq \frac{1}{2} \quad \mbox{and} \quad  \frac{4 \ell -1}{4\ell-2} \leq y\leq \frac{1}{2} \quad .
$$
Analogously, the neutral region $\mathcal{C}_{inv}$ corresponds to 
$$
\frac{1}{2} \leq x \leq \frac{3}{4} \quad \mbox{and} \quad \frac{1}{2}\leq y\leq 1  \quad .
$$
As expected, the number of time-steps needed for a trajectory to be captured 
by an attractor sensibly increases as the parameter $\ell$ approaches the equilibrium value $\ell=0.25$.
Figure \ref{bas} shows the variation of the size of the basins of attraction for increasing values of $\ell$ 
in the range $[0,\frac{1}{4}]$, where $\ell=\frac{1}{4}$ corresponds to the condition of ``equilibrium'' 
and $\ell=0$ corresponds to the most dissipative dynamics. In the latter case, apart from the ``neutral'' region $\mathcal{C}_{inv}$, which does not depend on $\ell$, the rest of the phase space constitutes the basin of attraction of the fixed point $P_\mathcal{D}$.
%In figure \ref{raz}, we show some specific examples corresponding to rational value of the parameter $\ell$: $\ell=\frac{1}{4}-\frac{1}{n}$.\\

\section{Phase Space contraction in the steady state}
As previously observed, the dynamics is neither ergodic nor chaotic and the Lyapunov exponents are either negative or vanish almost everywhere, hence fluctuations of the phase space contraction rate are prevented in the steady states.
As the phase space decomposes in several invariant sets, one may consider different invariant measures supported on the whole phase space or restricted to the single invariant sets. We refer to $\langle \Lambda \rangle_\mathcal{U}$ as the full phase space 
ensemble average which in general does not correspond to the time average computed along a single trajectory.
Indeed, in the steady state, starting from a uniform distribution of points, only certain values of the time-average phase space contraction rate are possible, each corresponding to a different attractor. We will denote each of the invariant sets by $\Omega_i$, and its corresponding basin of attraction by $\Omega_{0,i}$, with index $i=1,2,...N_\ell$, where $N_\ell$ is the number of invariant sets at a given $\ell$.

For instance, taking $\ell \geq \frac{1}{8}$, we have $N_\ell=6$, and there are $N_\ell-1$ different values for $\Lambda$ in the steady state:

\bea
\Lambda_{i}(\ell)  =
\left\{\begin{array}{cc}
\Phi,&   \Omega_{1}(\ell)=P_\mathcal{D}\\
\Phi/2,&  \Omega_{2} (\ell)={\mathcal{CDCD}}\\
0,&  \Omega_{3}(\ell)= \mathcal{B}_{inv}\\
0, & \Omega_{4}(\ell)=\mathcal{C}_{inv}\\
-\Phi/2,& \Omega_{6}(\ell)=\mathcal{AB}\\
-\Phi, & \Omega_{5}(\ell)=P_\mathcal{A}
\end{array}
\right. \quad .
\eea
For $\ell<\frac{1}{8}$, instead, it holds $N_\ell=4$:
\bea
\Lambda_{i}(\ell)  =
\left\{\begin{array}{cc}
\Phi,&  \Omega_{1}(\ell)=P_\mathcal{D}\\
0,& \Omega_{2}(\ell)= \mathcal{B}_{inv}\\
0, &  \Omega_{3}(\ell)=\mathcal{C}_{inv}\\
-\Phi, &  \Omega_{4}(\ell)=P_\mathcal{A}\\
\end{array}
\right. \quad .
\eea

If we compute the steady state average $\langle\Lambda\rangle $ over the whole phase space  $\mathcal{U}$, from a uniform initial distribution on $\mathcal{U}$, we obviously have:
\be
\langle\Lambda \rangle_{\ell,\mathcal{U}} = \sum_{i=1}^{N_{\ell}} \Lambda_i (\ell)\mu_0 (\Omega_{0,i} (\ell)) \quad ,
\ee
where  $\mu_0 (\Omega_{0,i} (\ell))$ denotes the Lebesgue measure of the $i$-th basin of attraction, which depends solely on $\ell$.

In our case, taking into account that the basins of $P_\mathcal{A}$ and $\mathcal{AB}$ have $0$ Lebesgue measure and, moreover, considering that the invariant regions $\mathcal{B}_{inv}$ and $\mathcal{C}_{inv}$ correspond to $\Lambda=0$, the following  holds:
$$
\langle \Lambda \rangle_{\ell,\mathcal{U}} =  - \ln(4\ell)  \left[  \mu_0(\Omega_{0,P_\mathcal{D}}(\ell))  +  \frac{1}{2} \mu_0(\Omega_{0,\mathcal{CDCD}})(\ell) \right] \quad ,
$$
because $\Phi = - \ln (\left|J_D\right|) = - \ln (4\ell)$.

\section{Transient and asymptotic FR }

In some sense, we may consider our map ``reversible'', because any trajectory of any number of steps $n$ can be associated with one trajectory producing the opposite phase space contraction rate, in $n$ steps. Furthermore, no pair of different trajectories needs to share the same reverse trajectory.
We obtained the initial conditions for the reverse trajectories by iterating the map for $n$ steps and then searching numerically for the initial conditions in $\mathcal{U}$ which lead to the  opposite $\Lambda$'s.
We numerically investigated, for different values of $\ell$, this ``reversibility" condition in order to test the transient FR.
Denoting by $P_n (\bar{\Lambda}_n = A) $ the probability (according to the uniform distribution) to get, over a $n$ step trajectory, the average phase space contraction rate $\bar{\Lambda}_n = A$ , the transient FR can be written in the following form:
\be
\frac{P_n (\bar{\Lambda}_n = A)}{P_n (\bar{\Lambda}_n = -A)} = e^{n \cdot A} \label{eqfrz} \quad .
\ee

\begin{figure}[ht!]
\center
\includegraphics[width=0.7\textwidth]{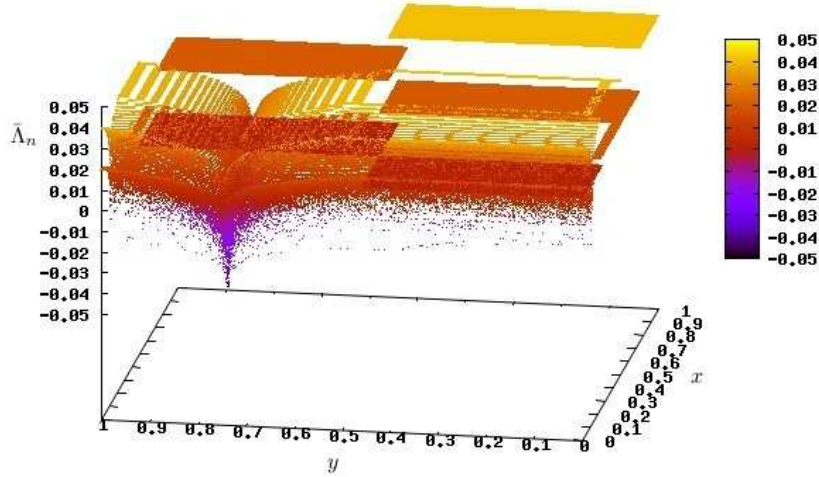}
\caption{Values of the phase space contraction time-average $\bar{\Lambda}_n$ for $\ell=0.24$ over trajectory segments of 250 steps. The initial conditions are $10^6$ uniformly randomly distributed points in the phase space.}\label{3D}
\end{figure}

If the transient FR holds, Eq.(\ref{eqfrz}) is valid for all $n\geq 0$. What happens, instead, in our model is that, for all values of $\ell$, 
Eq.(\ref{eqfrz}) is fulfilled for $n=1$ but not 
for $n=2$. On the other hand, increasing $n$, Eq.(\ref{eqfrz}) appears to be gradually restored. For small values of $\ell$, when the dynamics is highly dissipative, Eq.(\ref{eqfrz}) takes very few time-steps to be relatively accurately verified, while increasing $\ell$ towards $\frac{1}{4}$, it requires larger $n$ to reach the same accuracy.
\begin{figure}[htbp!]
\center
\includegraphics[scale=0.3]{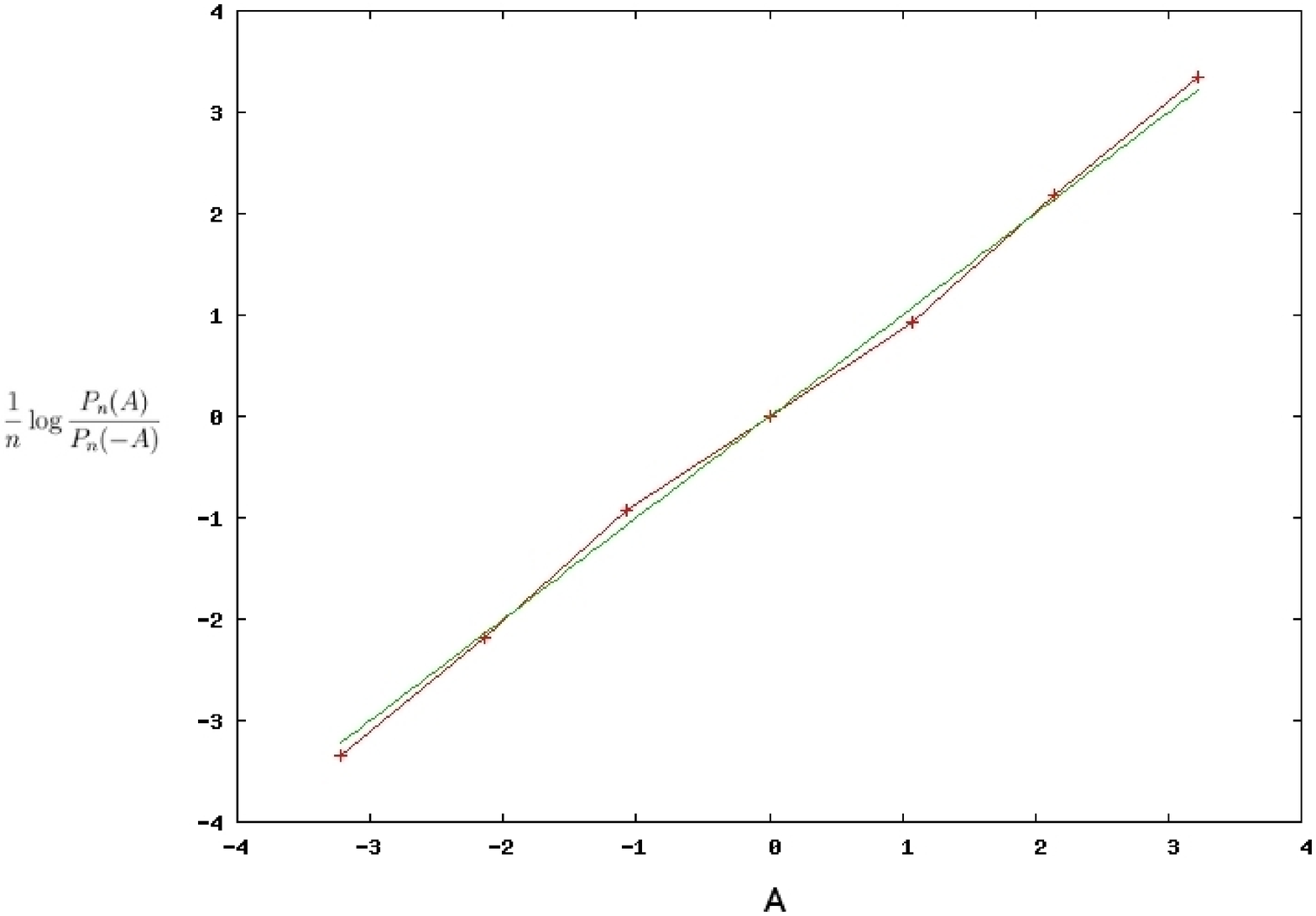} \quad
\includegraphics[scale=0.3]{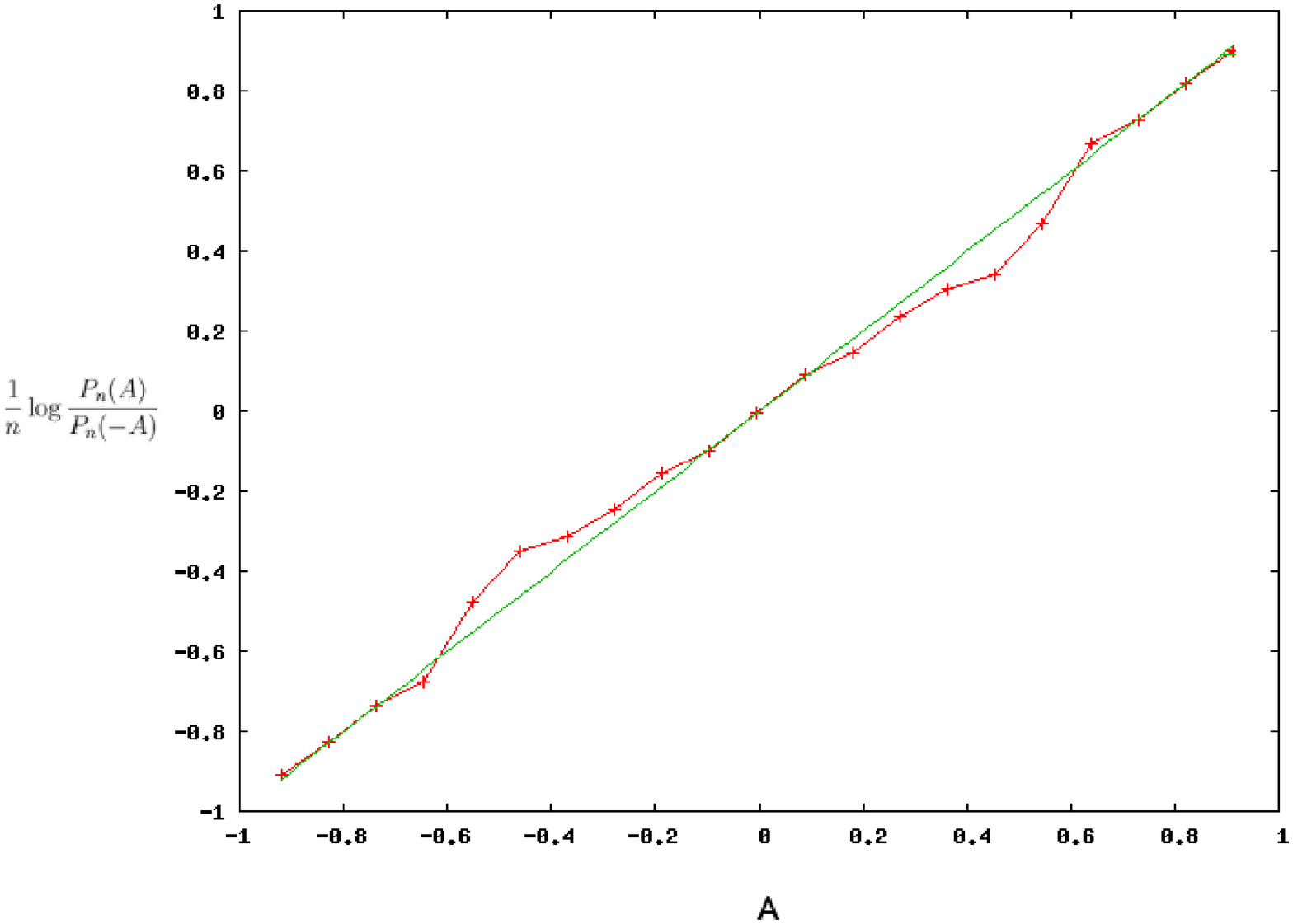}

\includegraphics[scale=0.3]{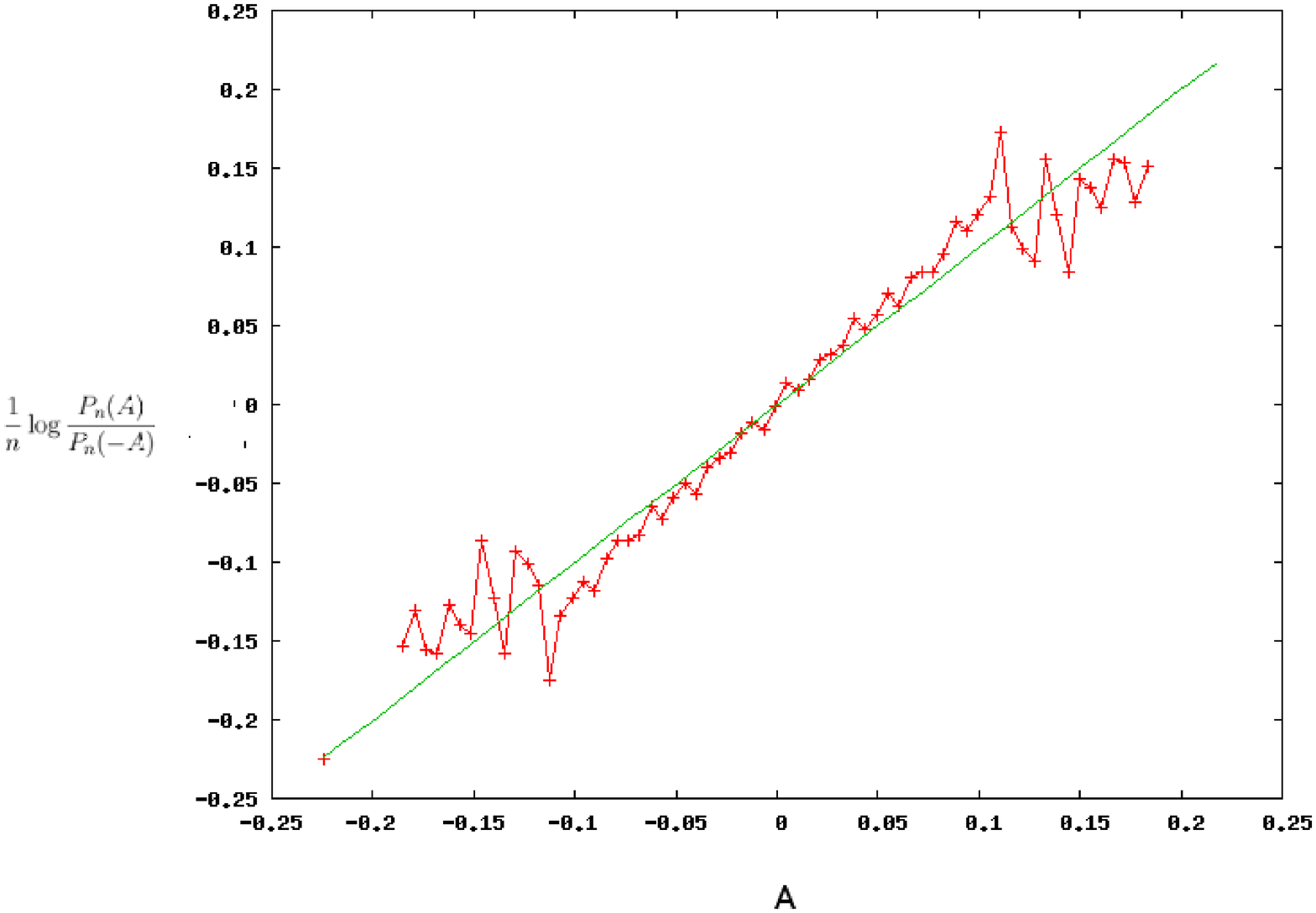} \quad
\includegraphics[scale=0.3]{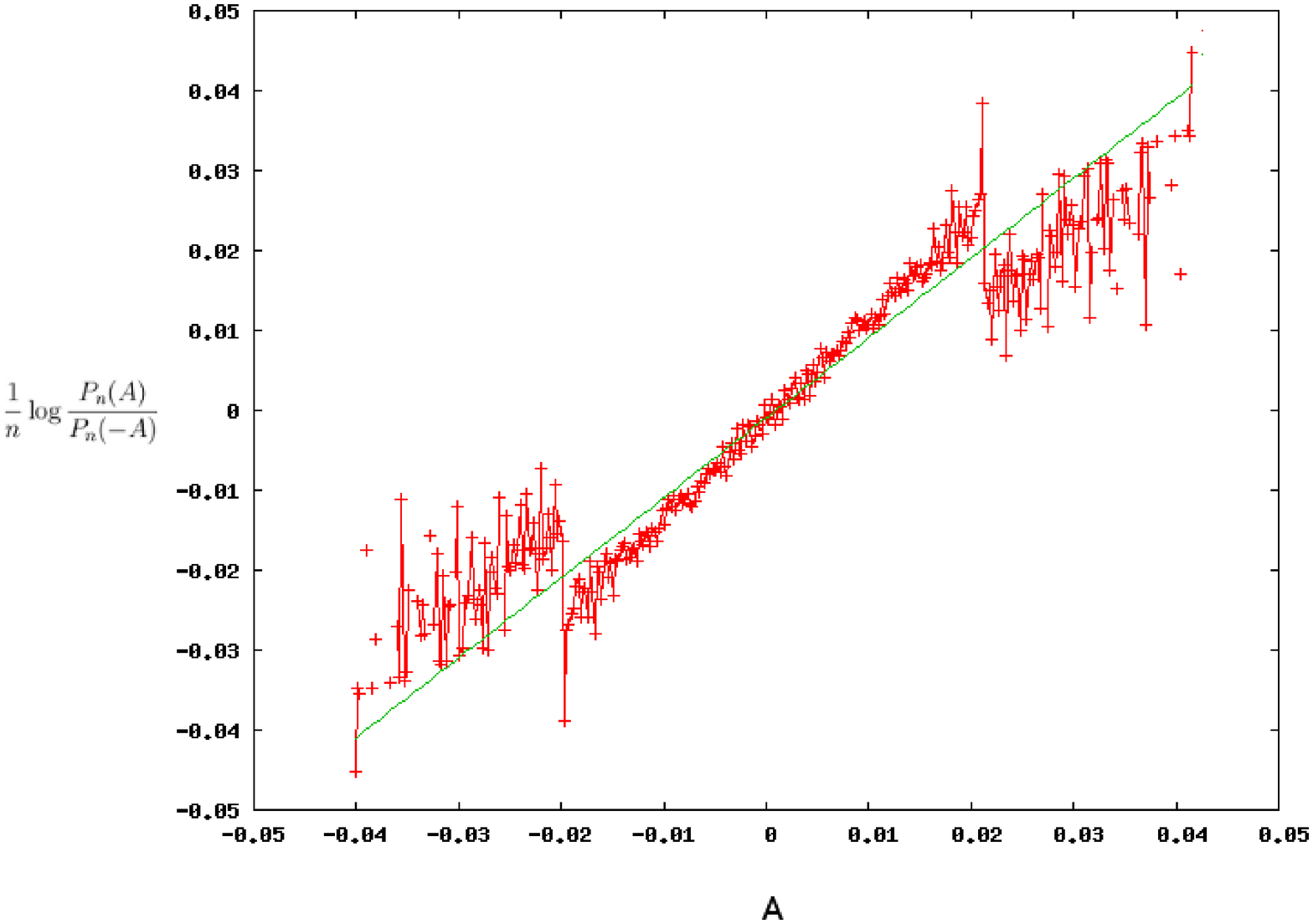}

\caption{FR restored in the long time limit for different $\ell$ values. From left to right, top to the bottom: $\ell=0.001$ and $n=3$, $\ell=0.1$ and $n=10$, $\ell=0.2$ and $n=40$; $\ell=0.24$ and $n=180$.}\label{figFR}
\end{figure}

As illustrated by Fig.\ref{3D}, starting from a point close to the repeller $P_\mathcal{A}$, the dynamics remains for a long time within the 
expanding region $A$, hence leading to the largest negative values of $\bar{\Lambda}_n$, for $n$ fixed.
\begin{figure}[htbp!]
\center
\includegraphics[scale=0.32]{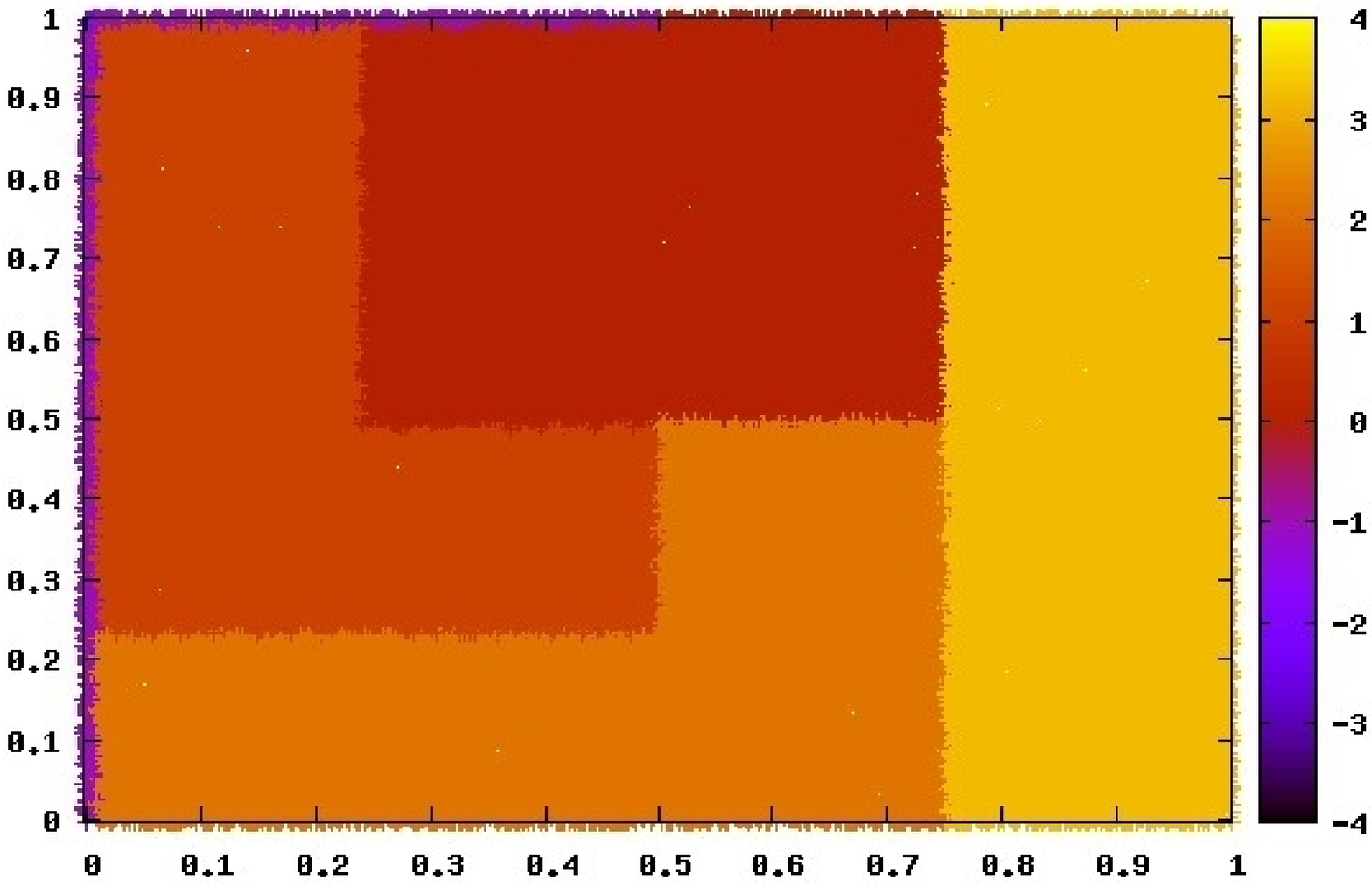}
\includegraphics[scale=0.32]{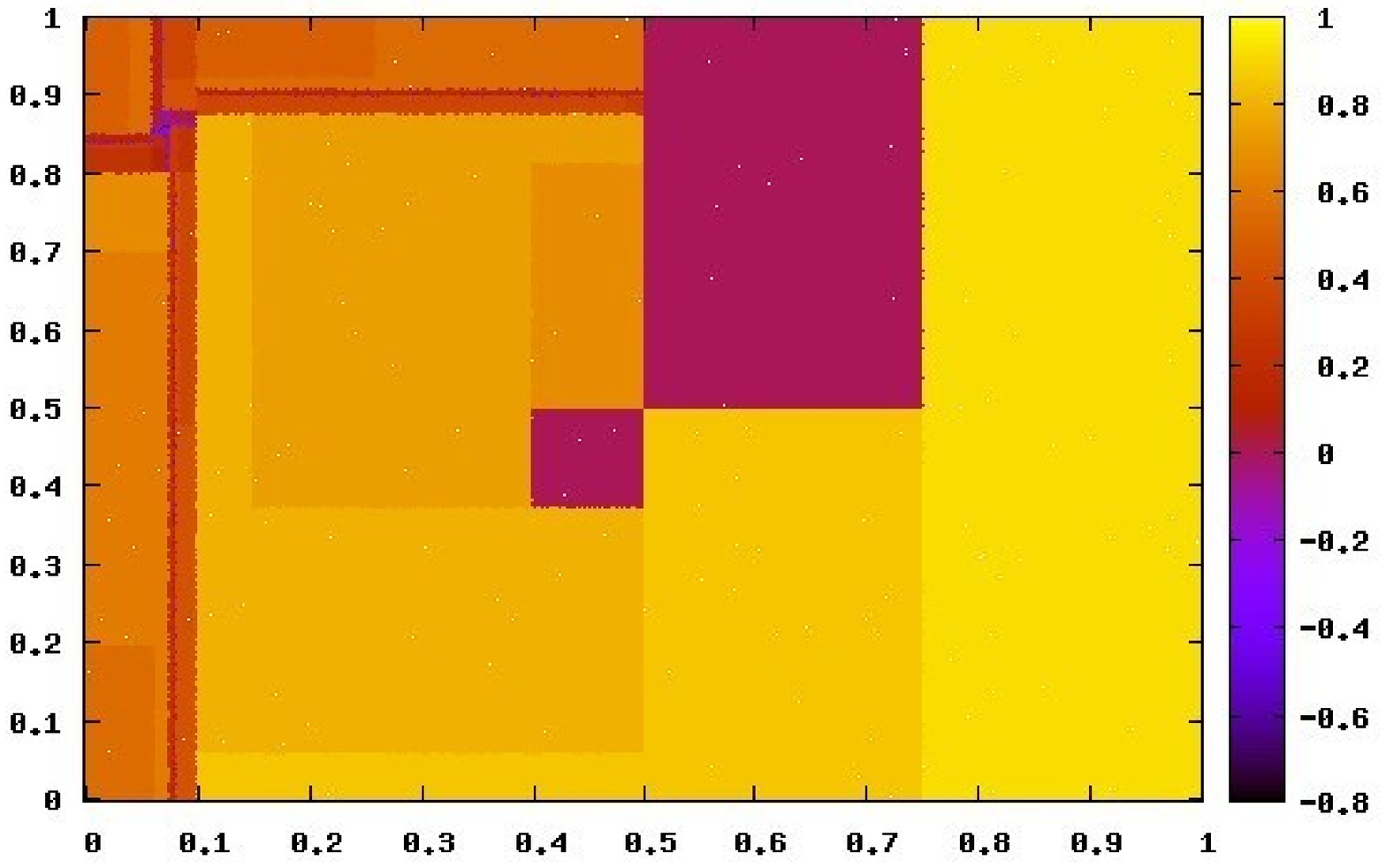}

\includegraphics[scale=0.32]{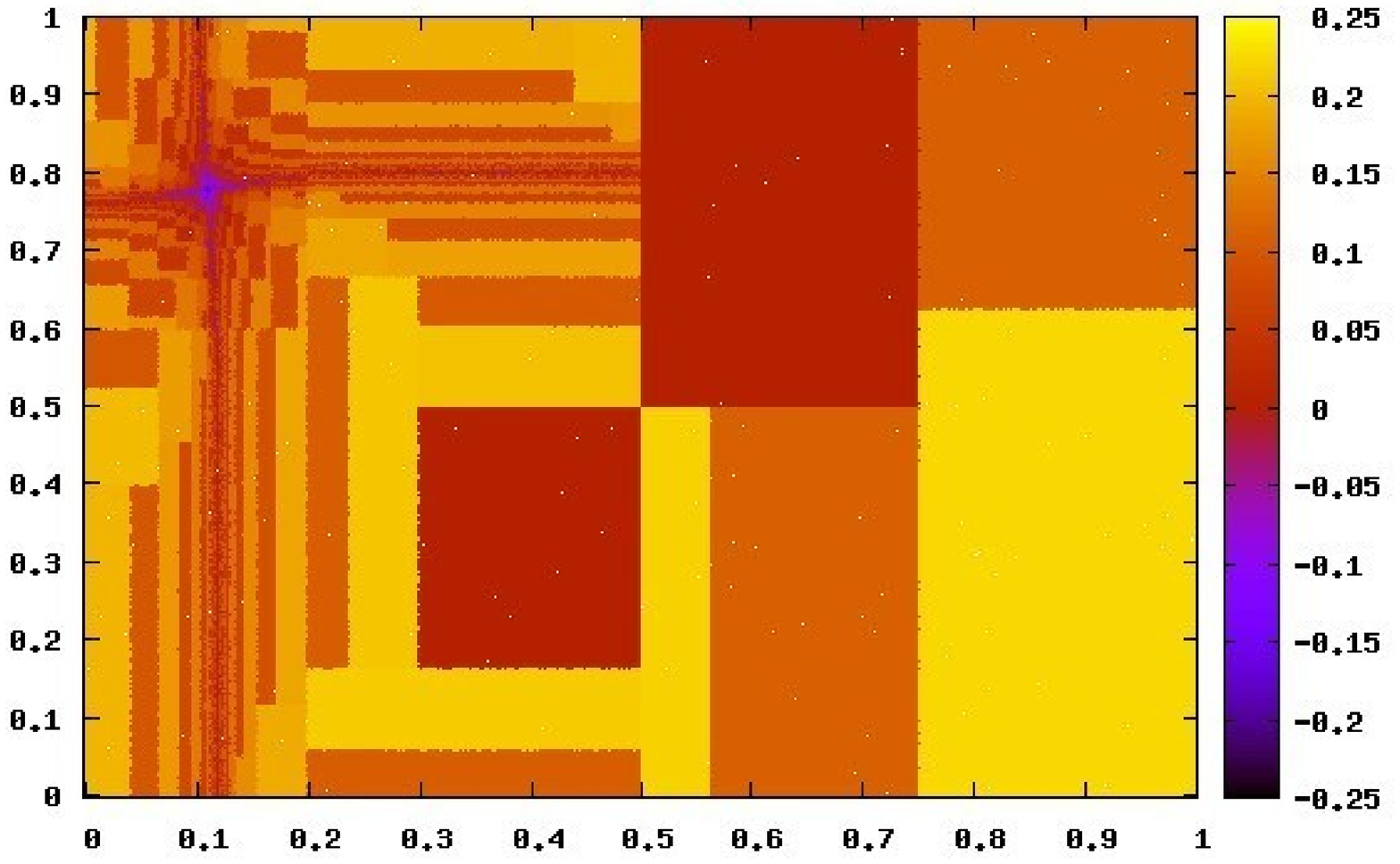}
\includegraphics[scale=0.32]{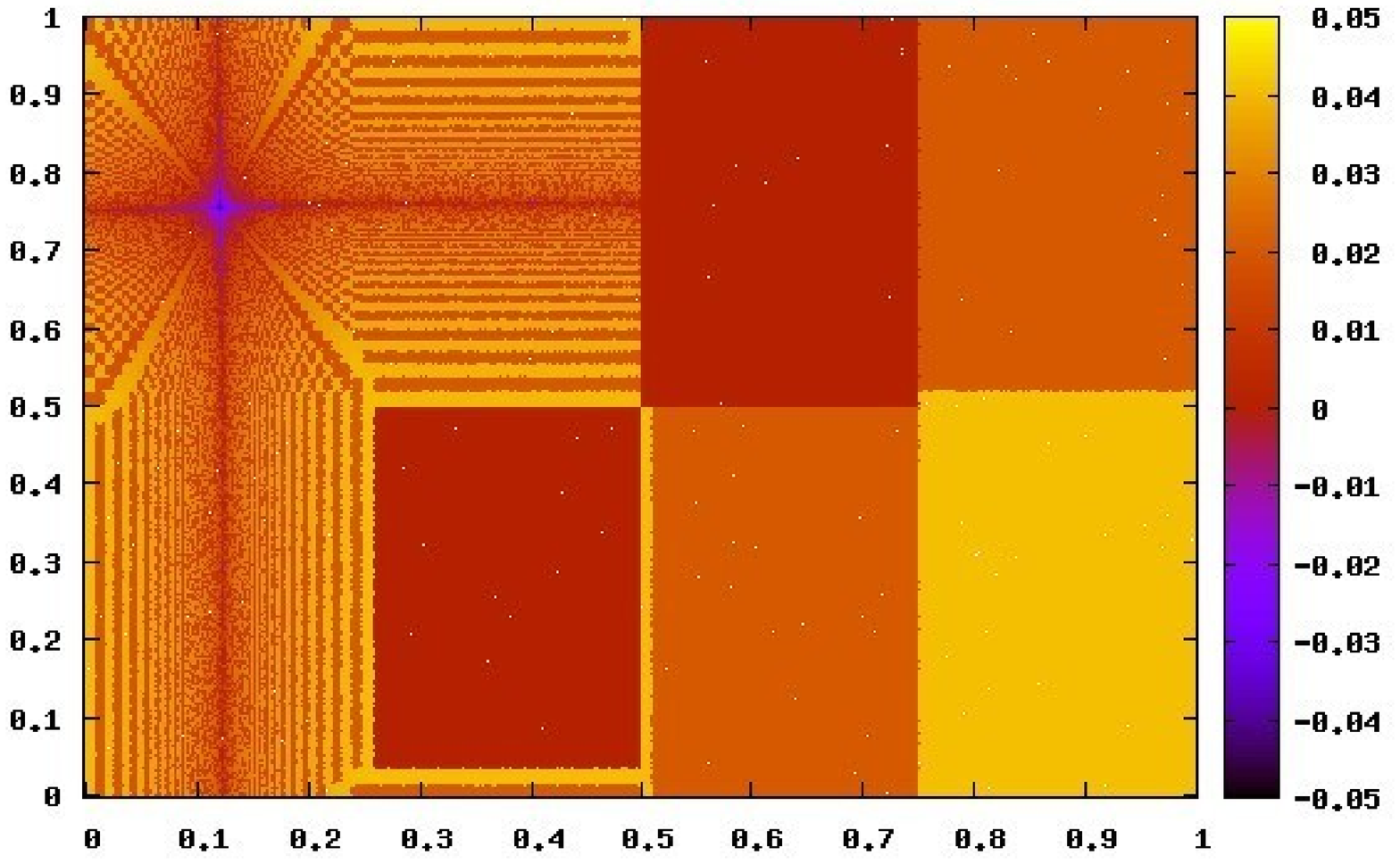}

\caption{Contour plot of the distribution of the average phase space contraction rate for increasing values of $\ell$. The dark black and purple colours around the repeller represent the locus of points which produce in the transient regime the most negative values. The length of the trajectory depends on the $\ell$ parameter and was chosen as the most suitable to validate the FR.
{\itshape From left to right, top to the bottom:} $\ell=0.001$ and $n=3$, $\ell=0.1$ and $n=15$, $\ell=0.2$ and $n=50$, $\ell=0.24$ and $n=250$.}
\label{fig11}
\end{figure}
Since the equilibrium dynamics is not ergodic, we should not expect the transient FR to be verified in general, because this ergodicity is required by the known derivations of the transient FR. Nevertheless our numerical results indicate that  Eq.(\ref{eqfrz}) is asymptotically fulfilled.
Fig.\ref{figFR}, \ref{fig11} and \ref{trfr} show the distribution of $\bar{\Lambda}_n$ on the map: on the left panel of Fig. \ref{figFR} and in Fig. \ref{trfr}, it can be observed that the values of the quantity
$$
\frac{1}{n} \log\left( \frac{P_n (A)}{P_n(-A)} \right)
$$
approach the line of slope 1 when $n$ increases, meaning that the distribution of $\Lambda$ tends to fulfill Eq. (\ref{eqfrz}) at large $n$.
\begin{figure}[htpb!]

\includegraphics[scale=0.20]{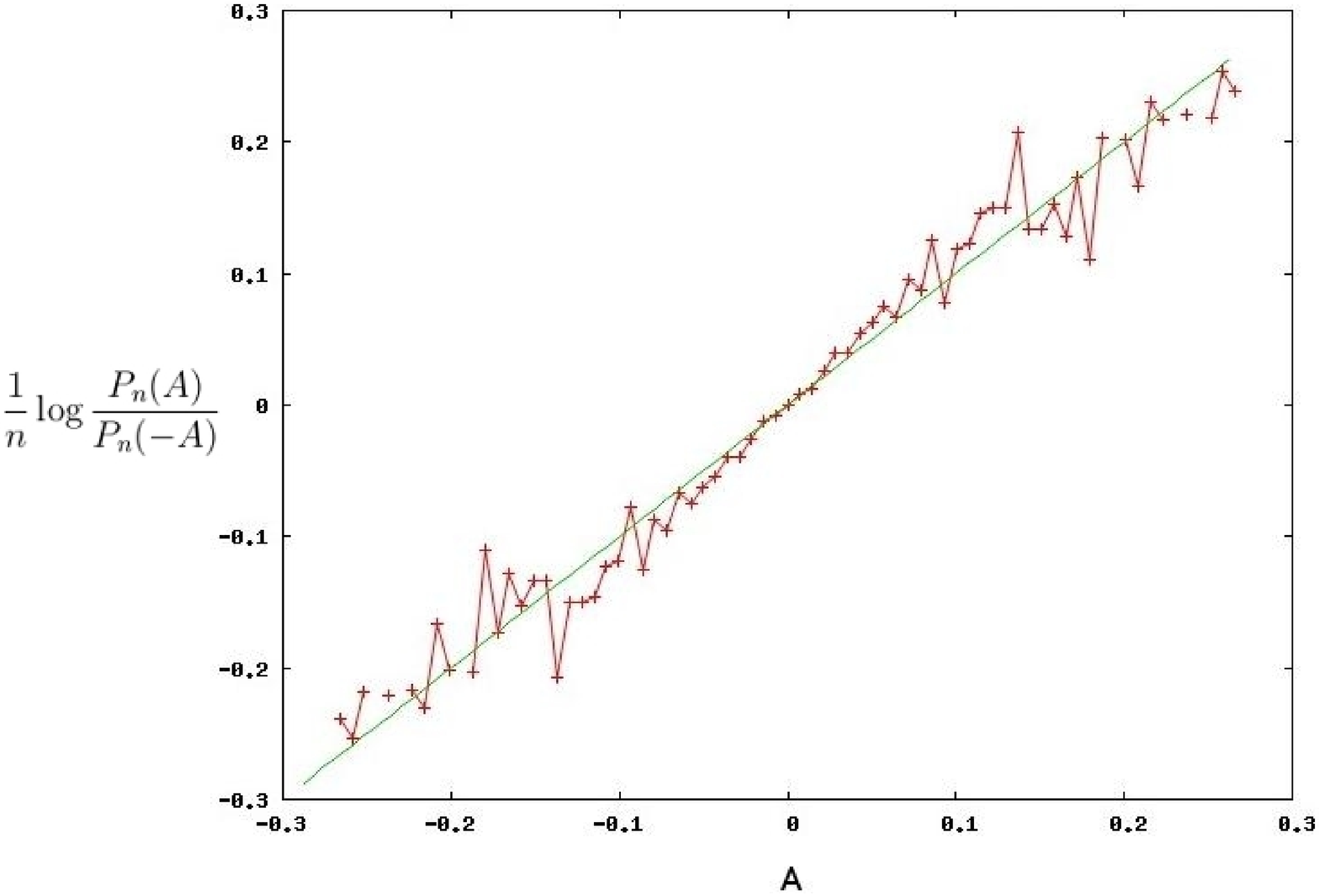}
\includegraphics[scale=0.21]{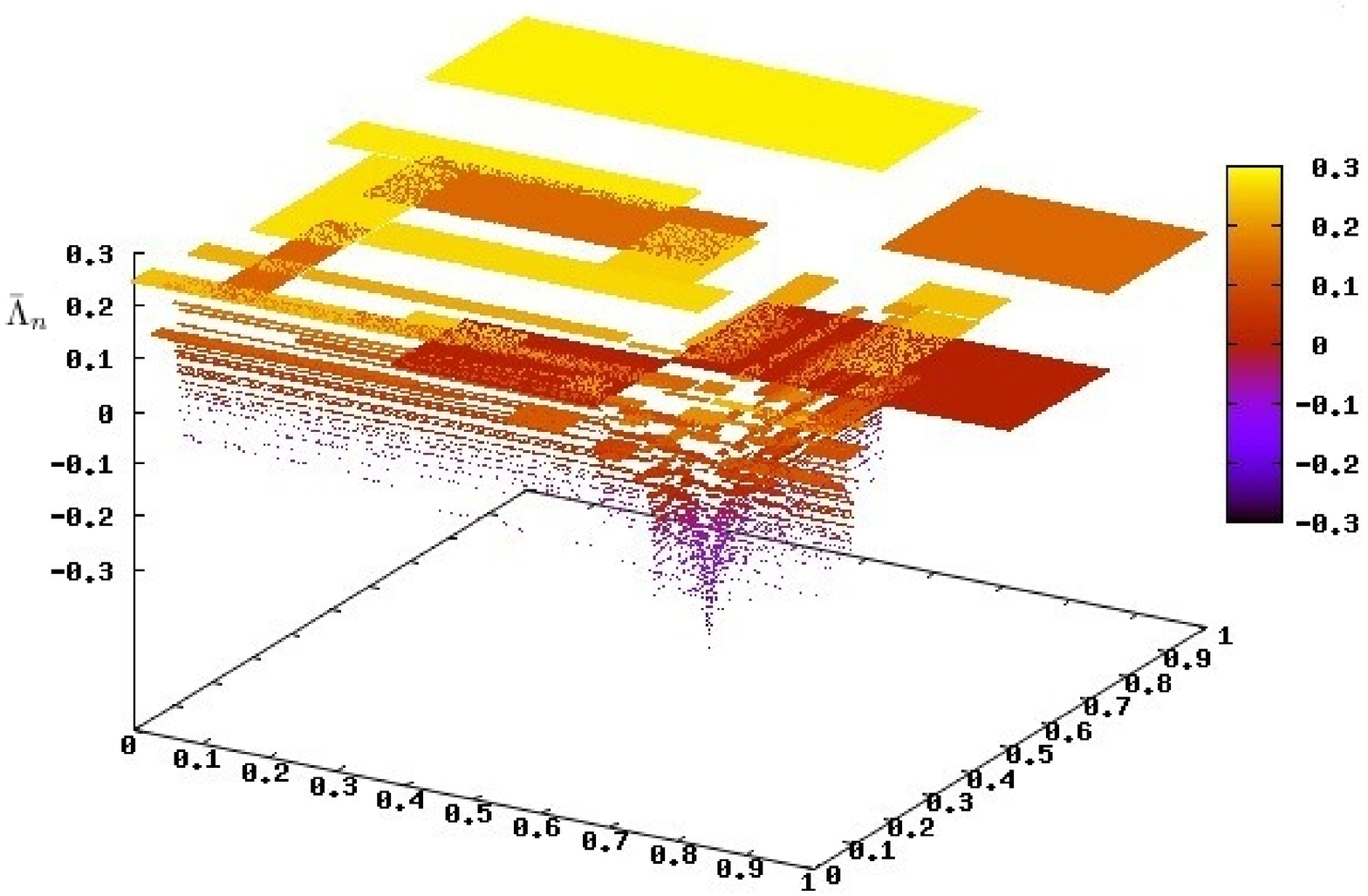}
\includegraphics[scale=0.21]{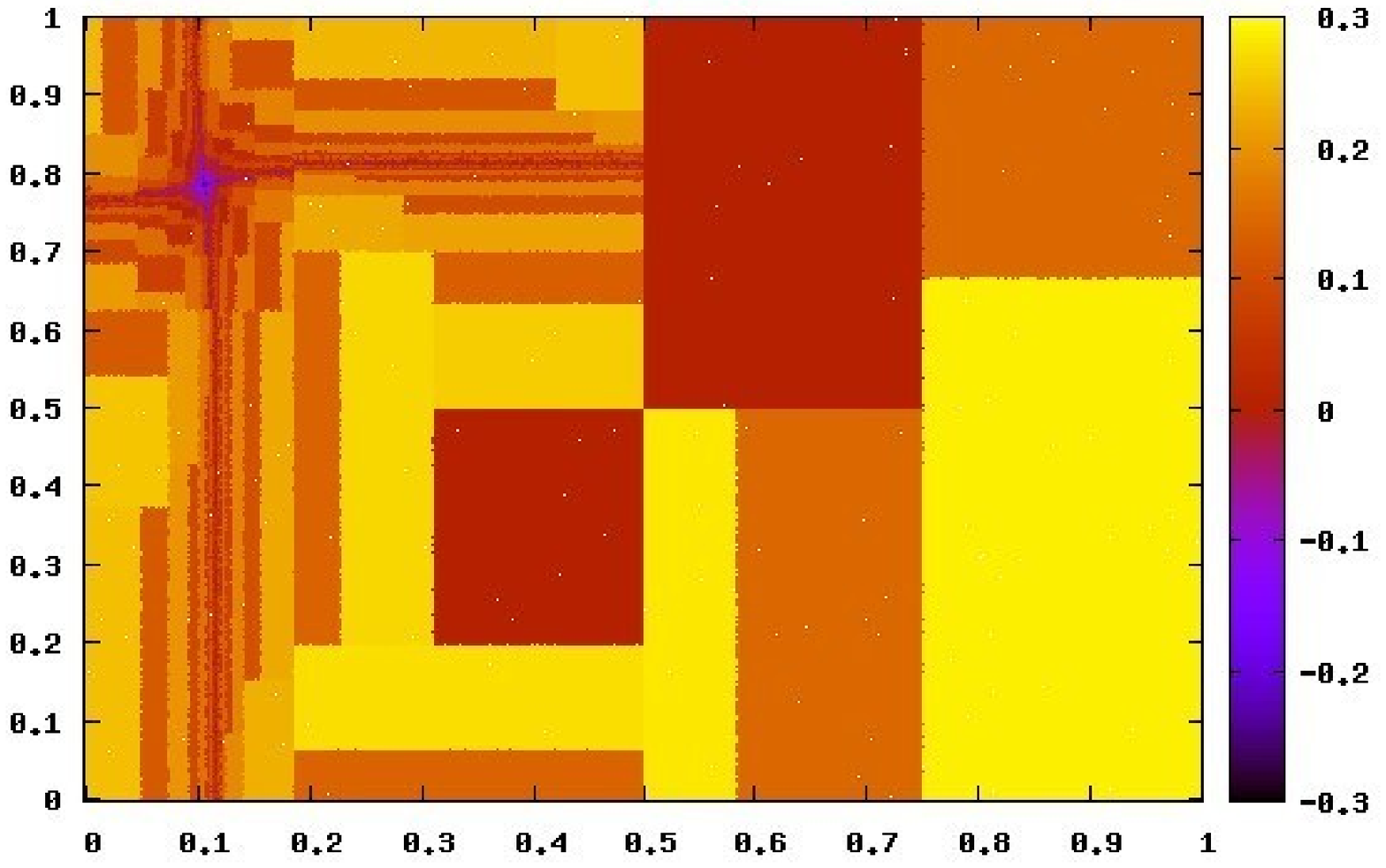}

\includegraphics[scale=0.20]{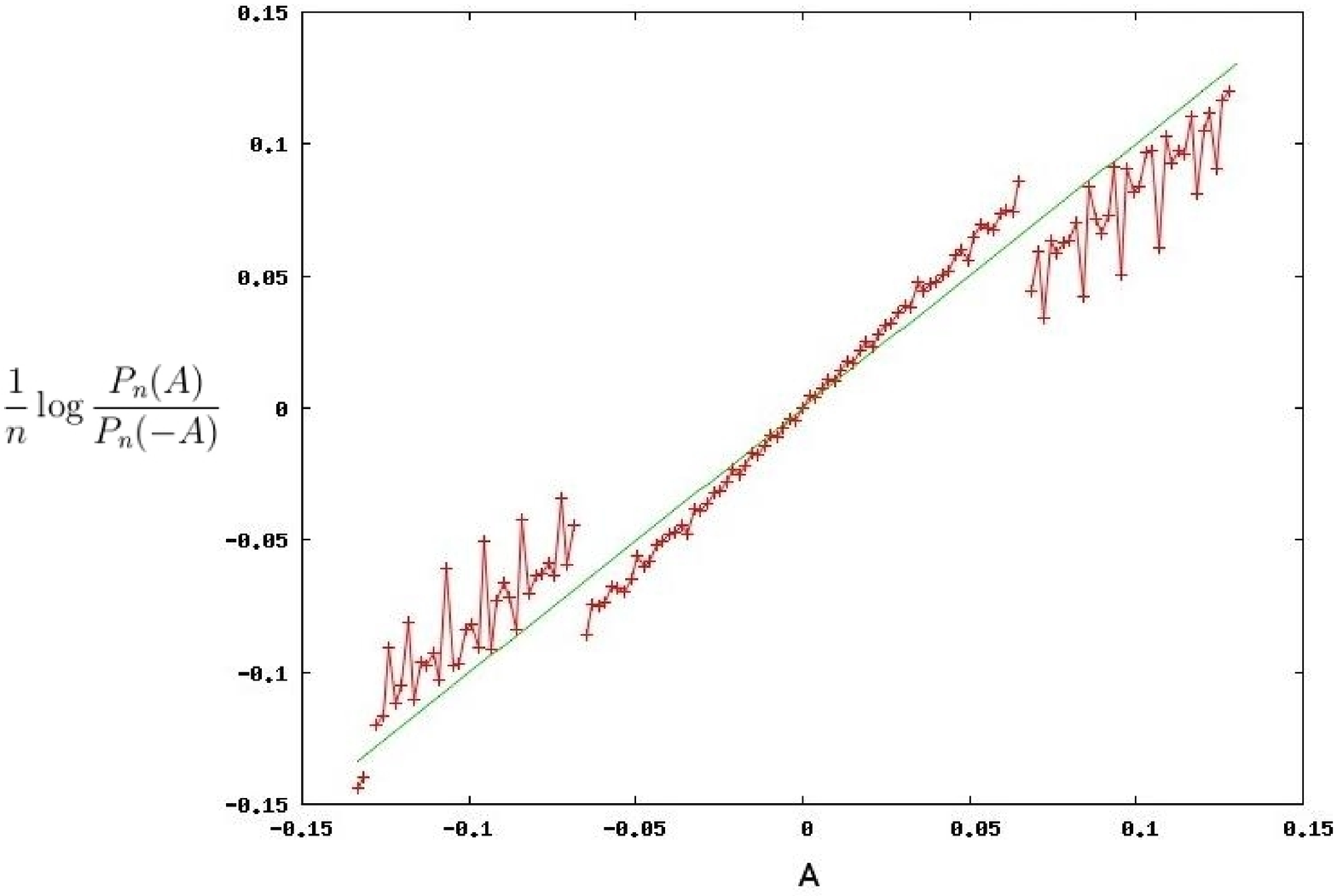}
\includegraphics[scale=0.21]{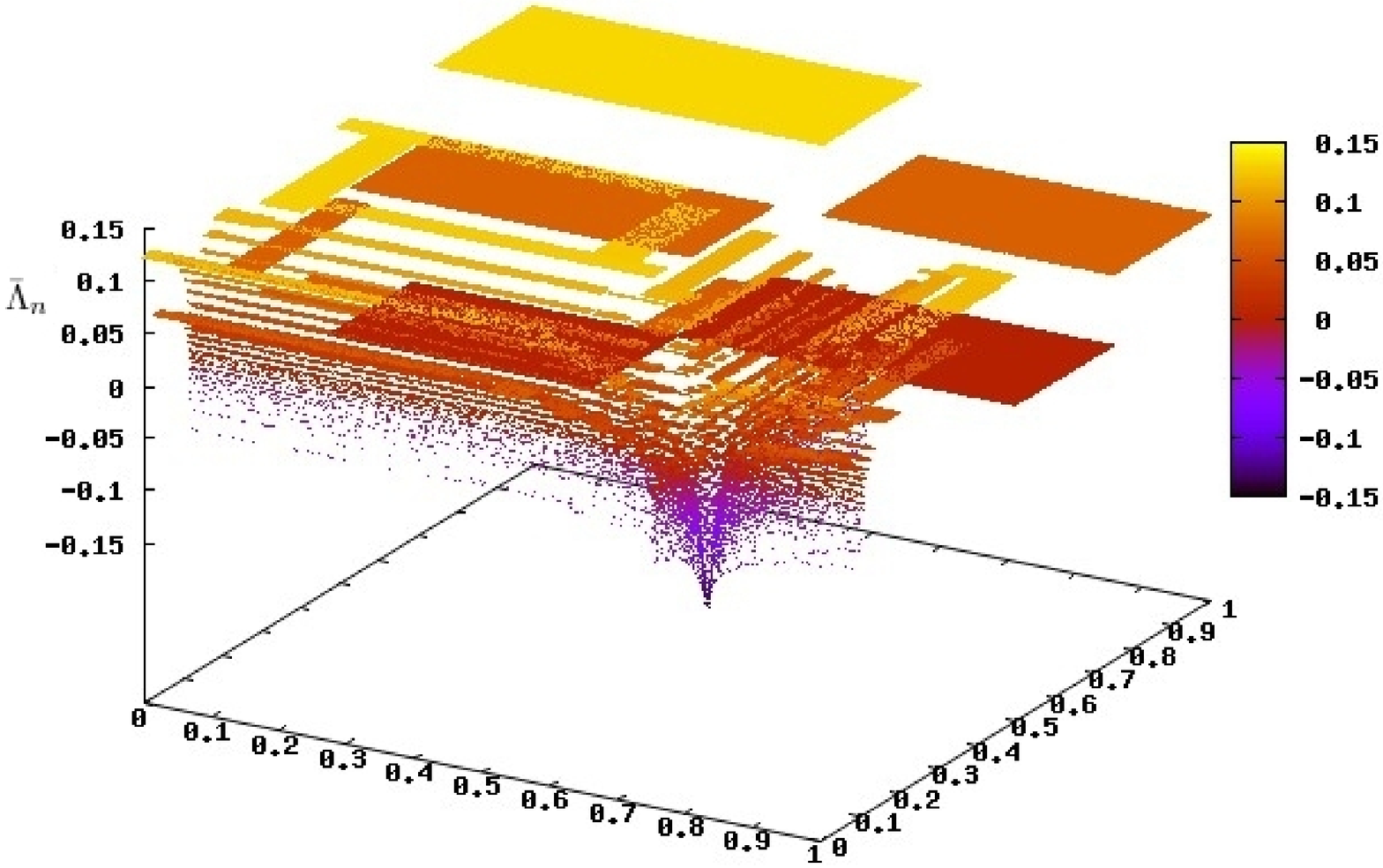}
\includegraphics[scale=0.21]{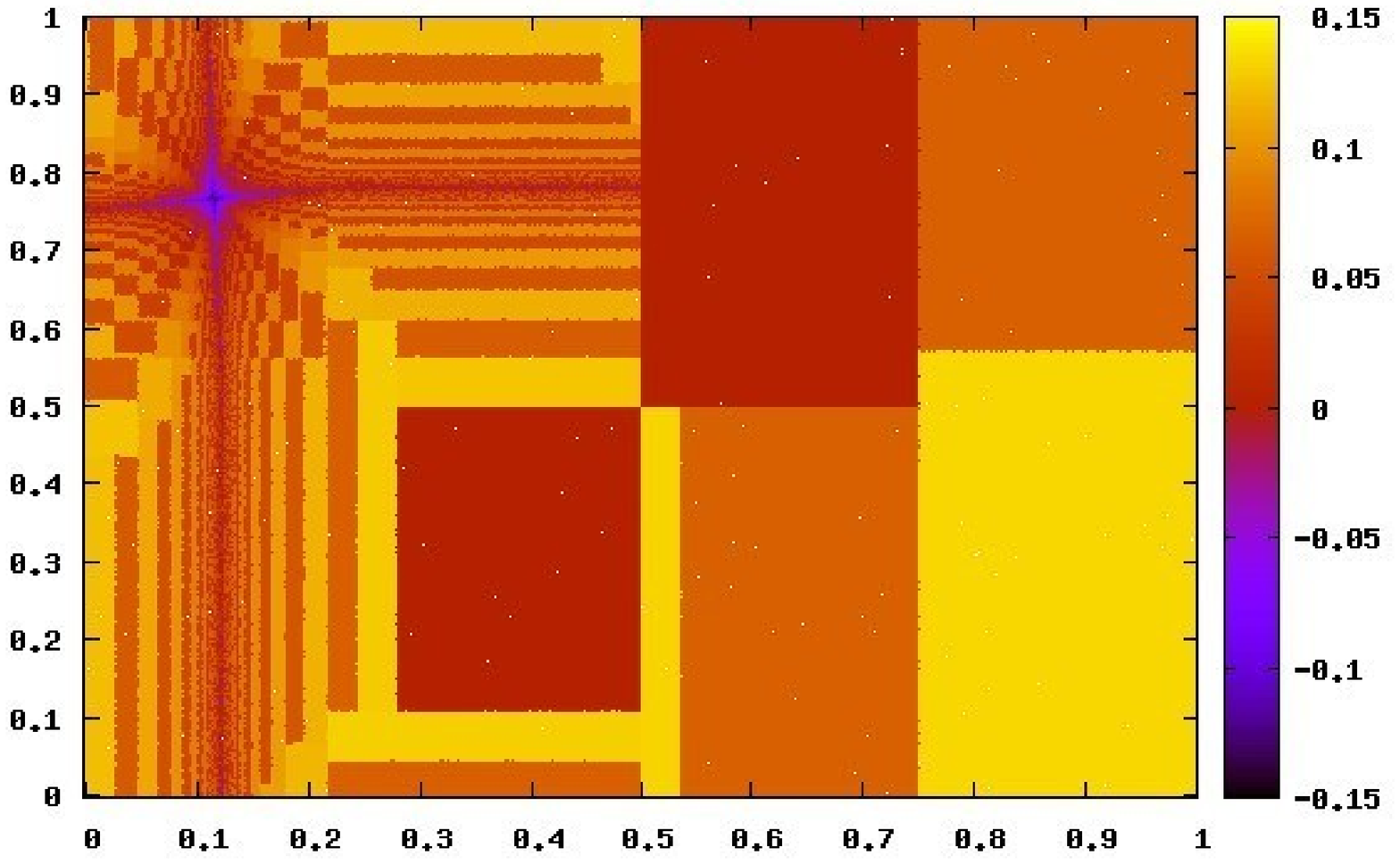}

\includegraphics[scale=0.20]{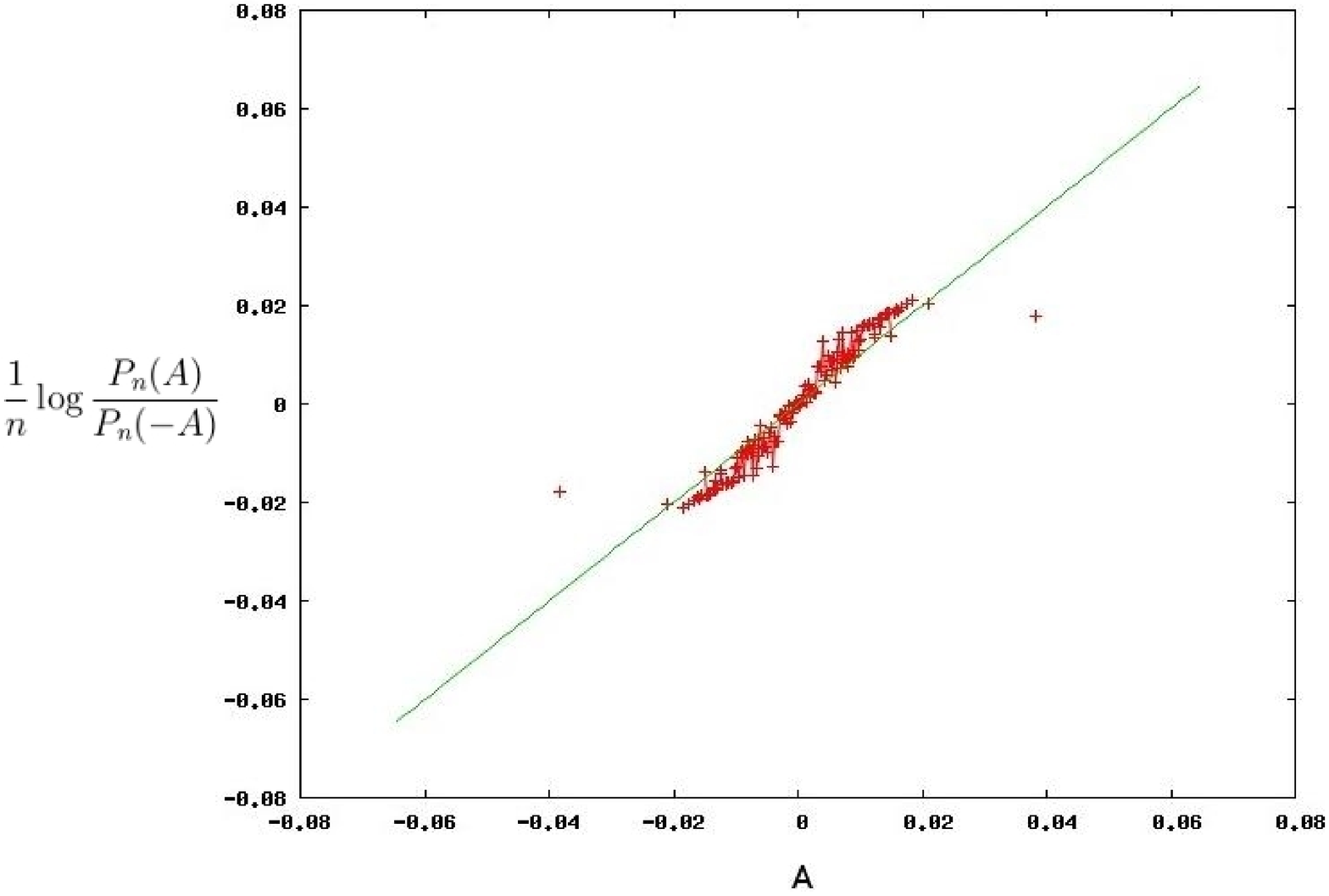}
\includegraphics[scale=0.21]{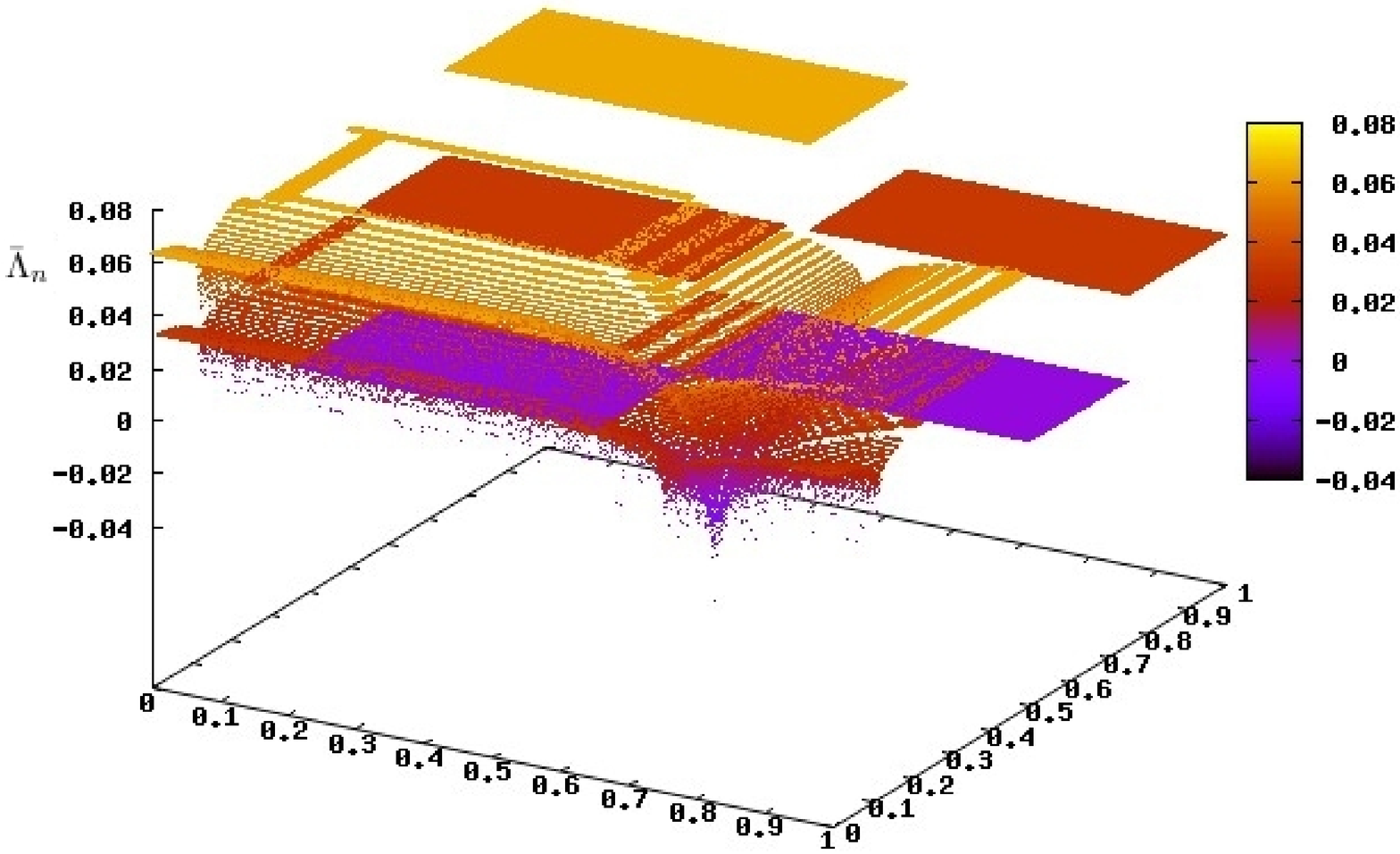}
\includegraphics[scale=0.21]{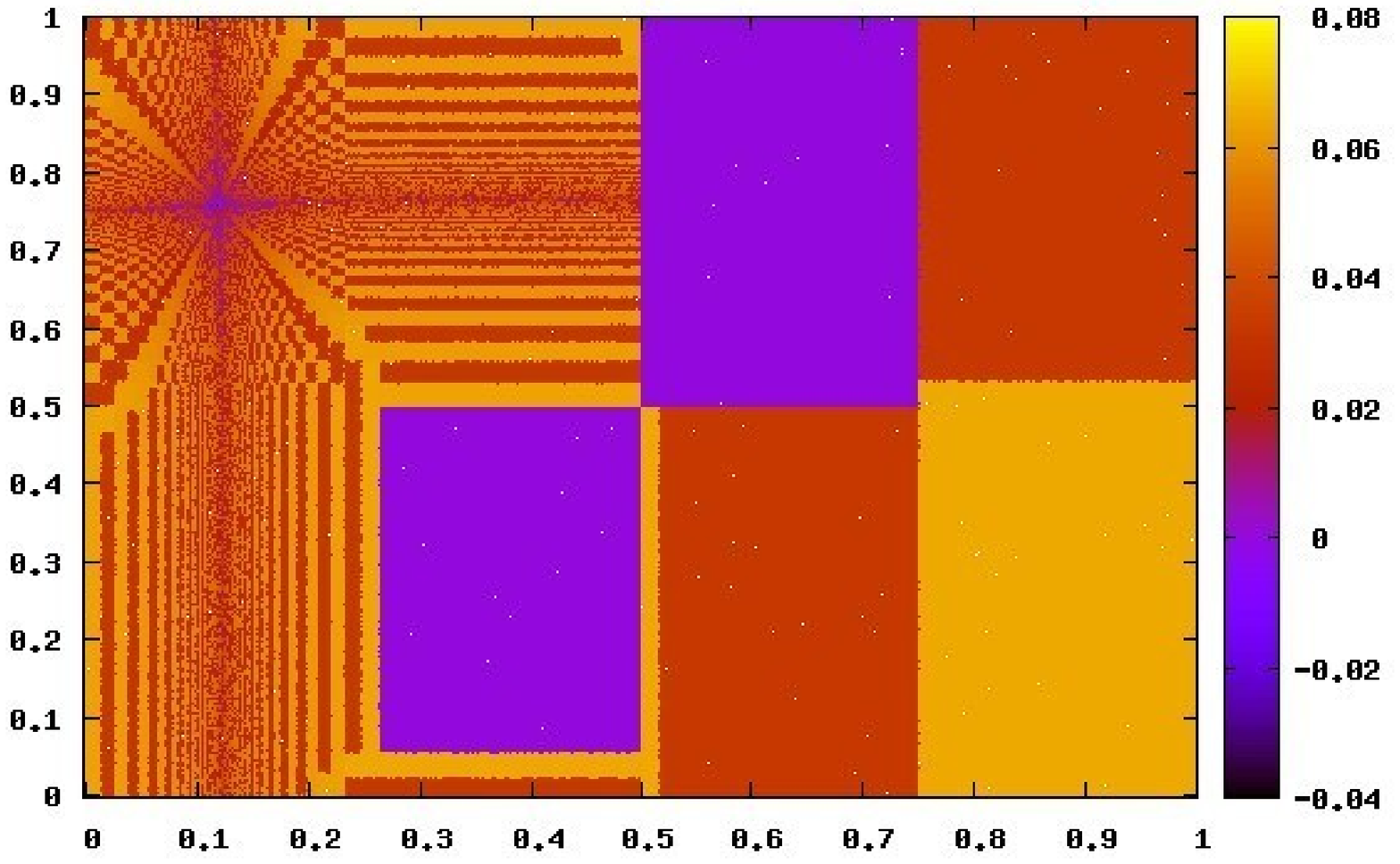}

\caption{On the left, approximate verification of the transient Fluctuation Relation. In the center, values of $\bar{\Lambda}_n$. On the right, the corresponding contour plot for (from top to the bottom): $\ell=\frac{1}{4}-\frac{1}{16}$ and $n=30$, $\ell=\frac{1}{4}-\frac{1}{32}$ and $n=70$, $\ell=\frac{1}{4}-\frac{1}{64}$ and $n=450$.}\label{trfr}
\end{figure}
Thus, what we observe is neither the standard transient nor the steady state FR, rather it is something that can be regarded as an 
asymptotic relation of a different kind, not yet reported in the literature. 

\begin{figure}[htbp!]
\includegraphics[scale = 0.18]{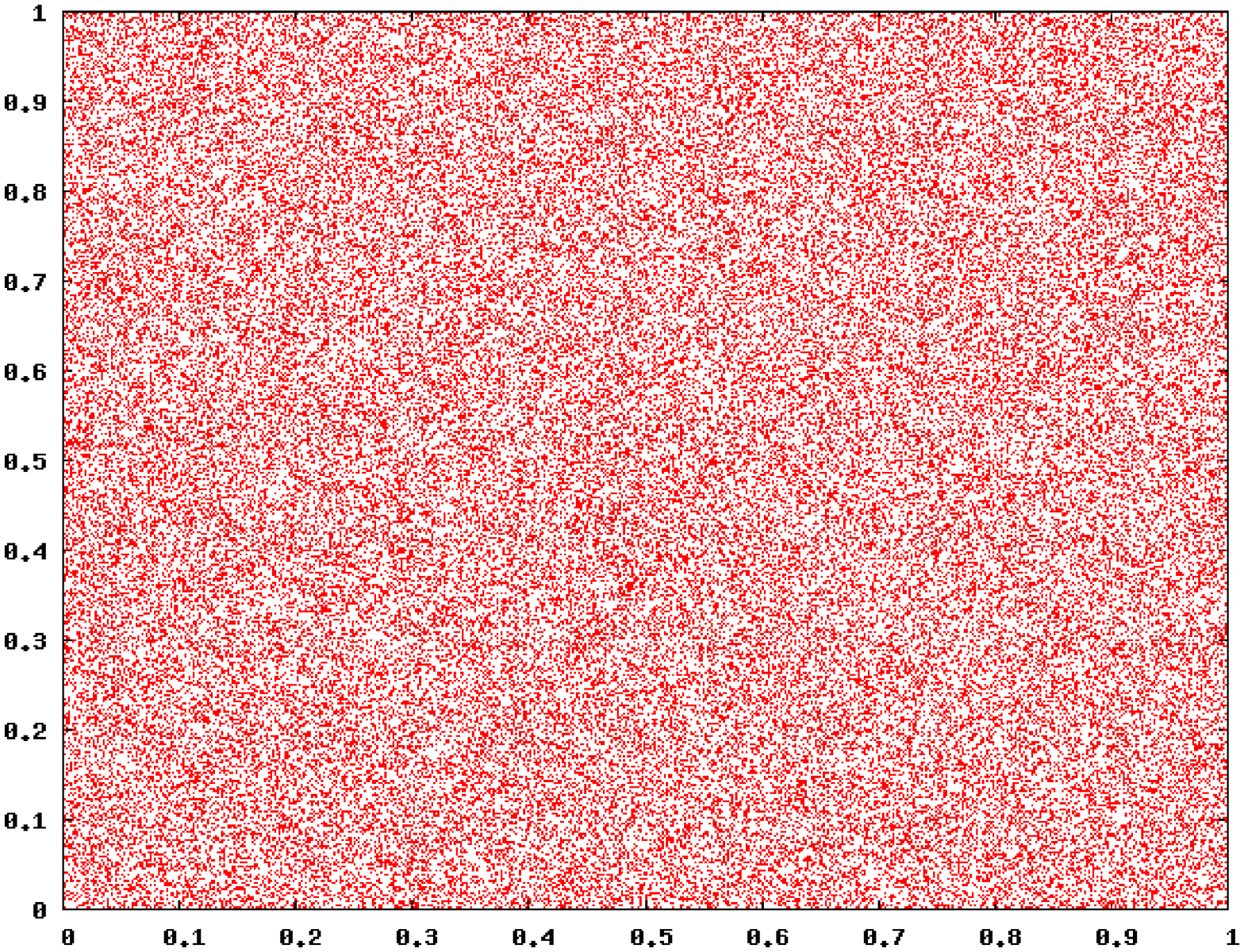}
\includegraphics[scale = 0.18]{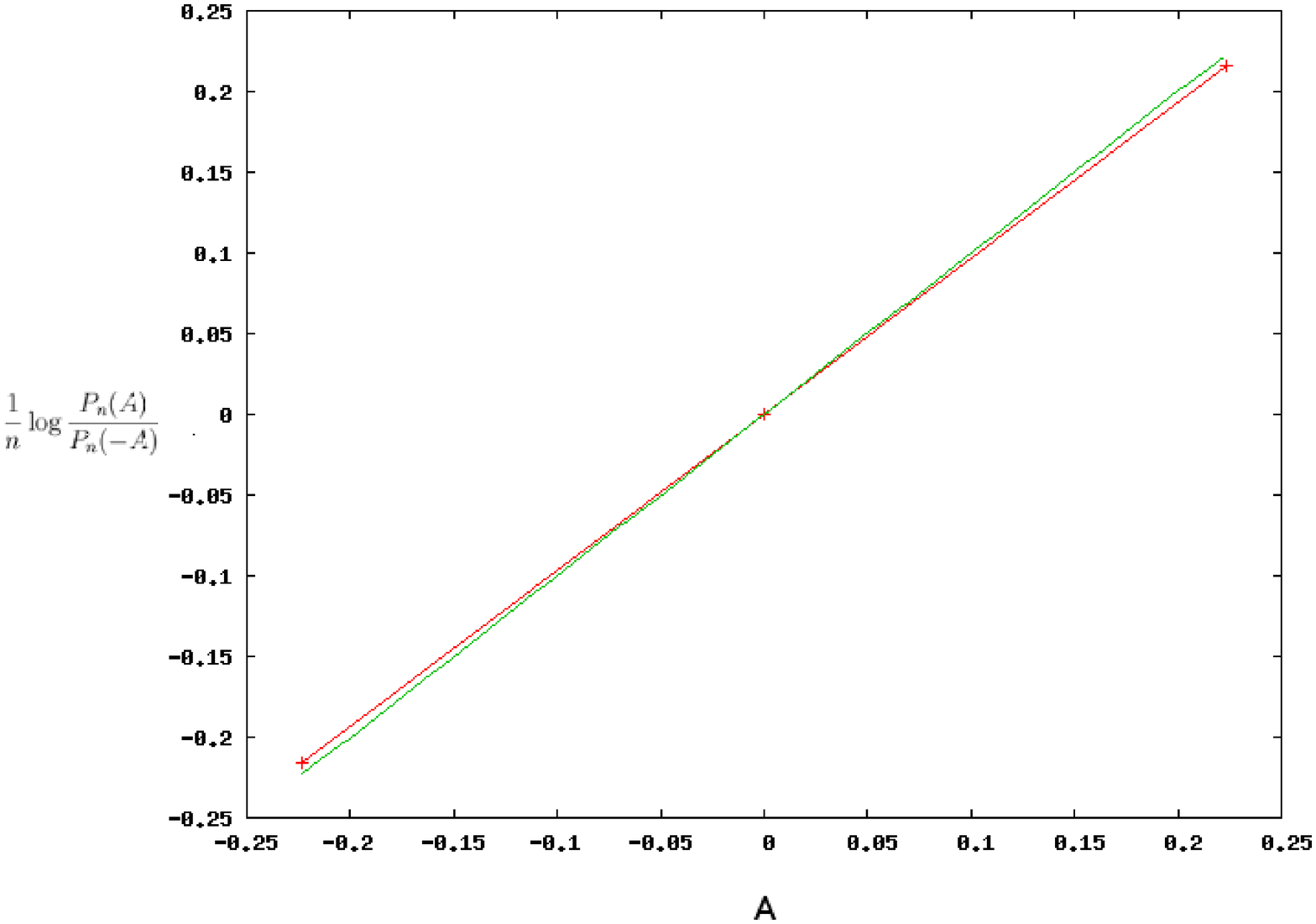}
\includegraphics[scale = 0.18]{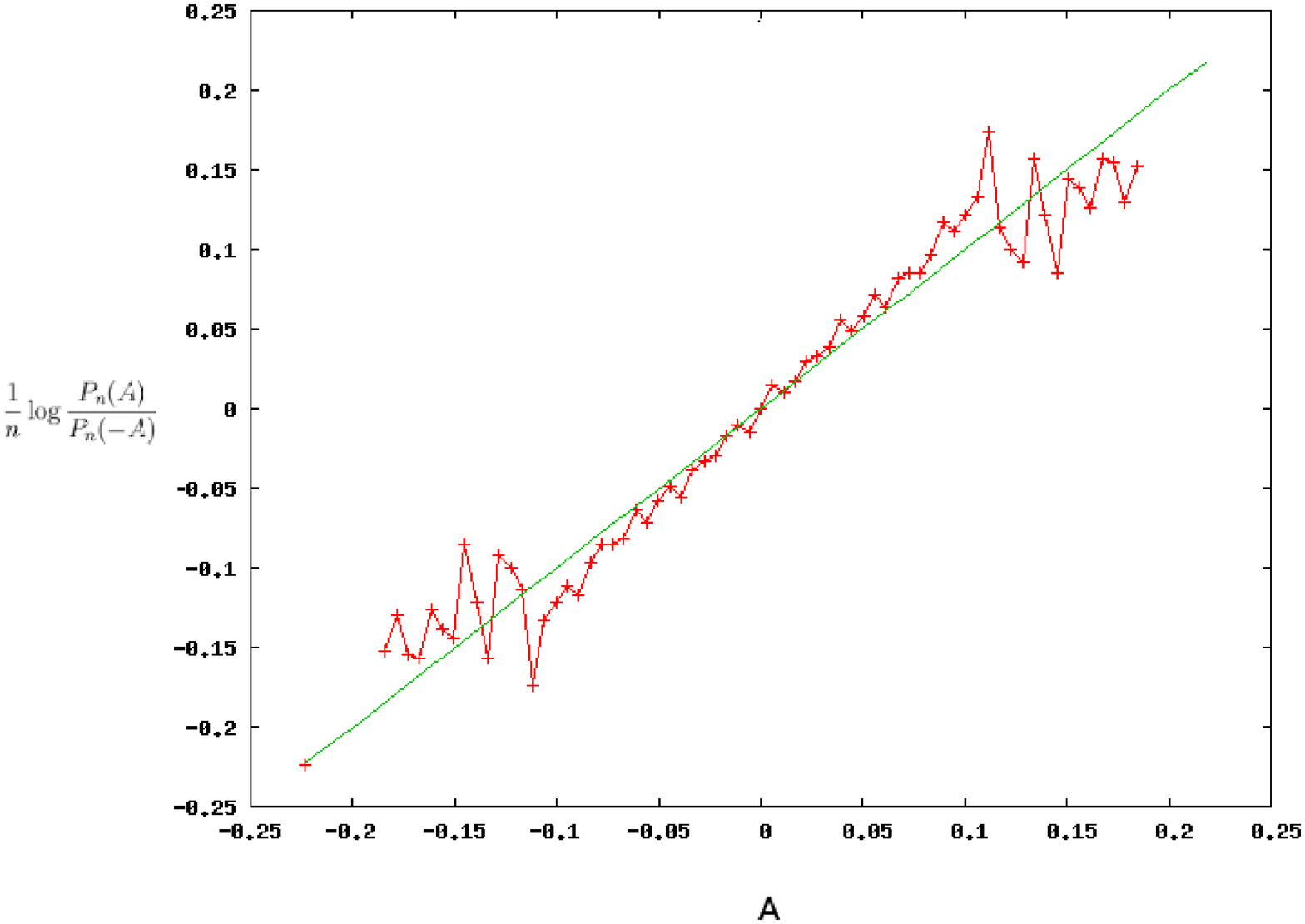}

\includegraphics[scale = 0.18]{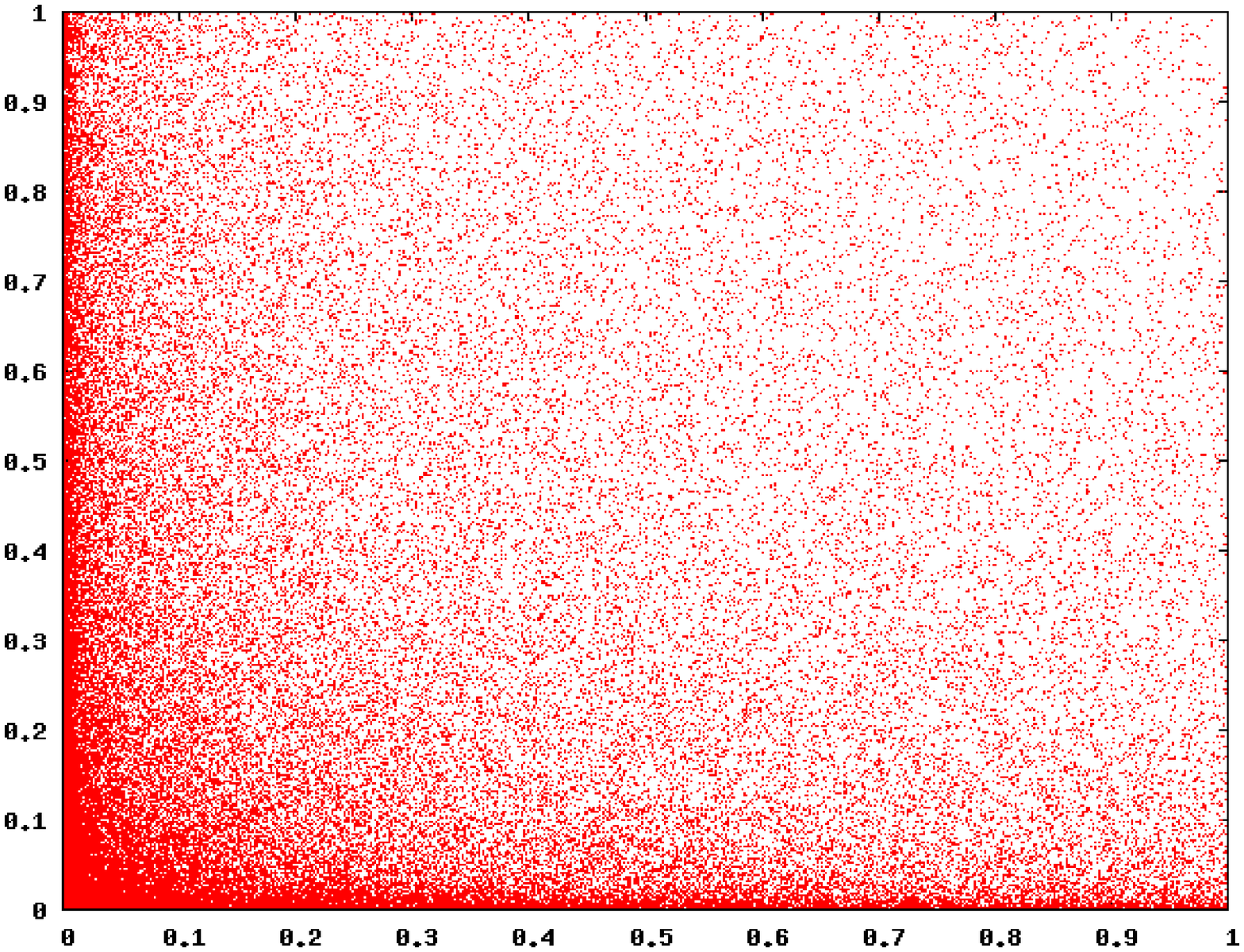}
\includegraphics[scale = 0.18]{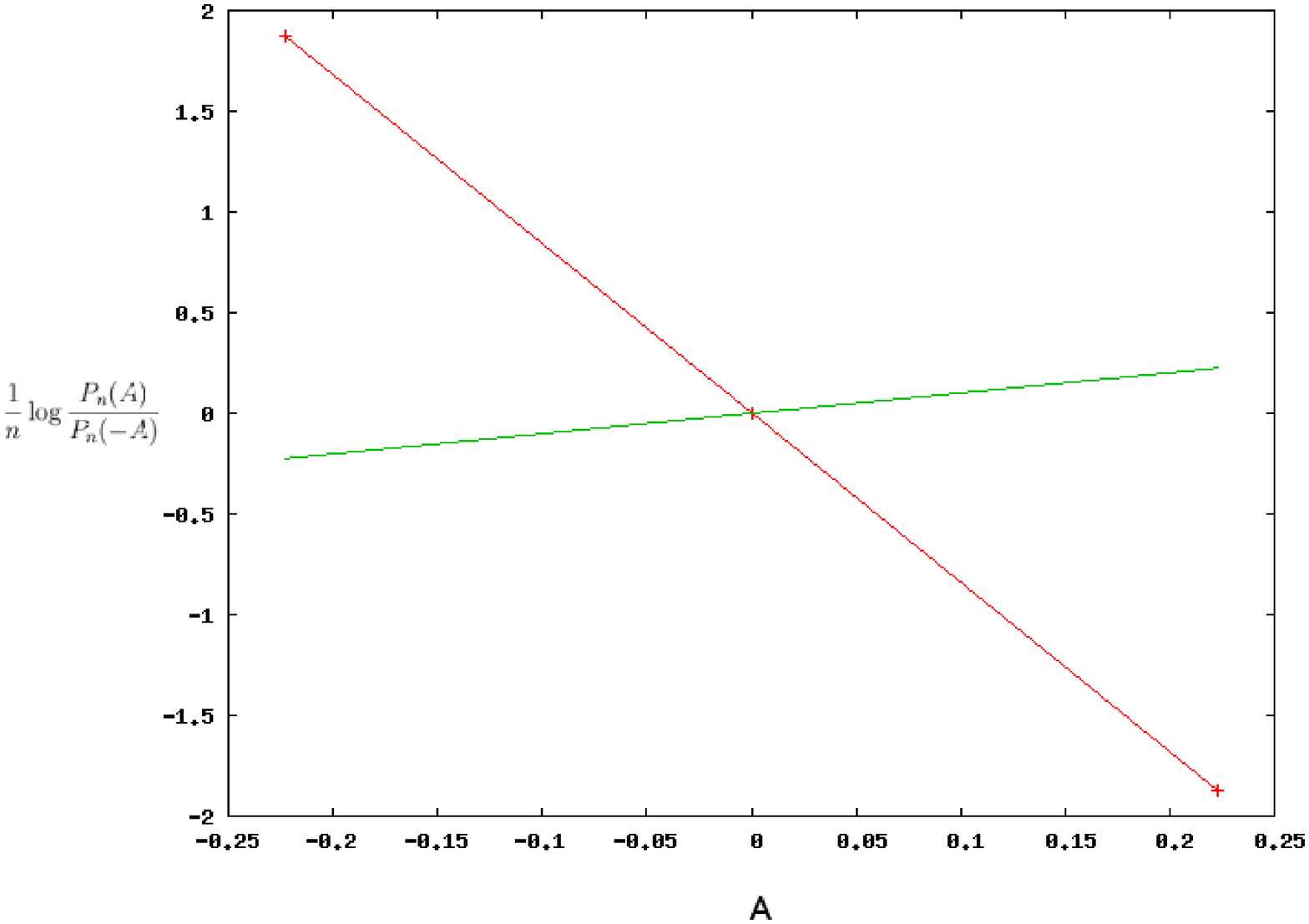}
\includegraphics[scale = 0.18]{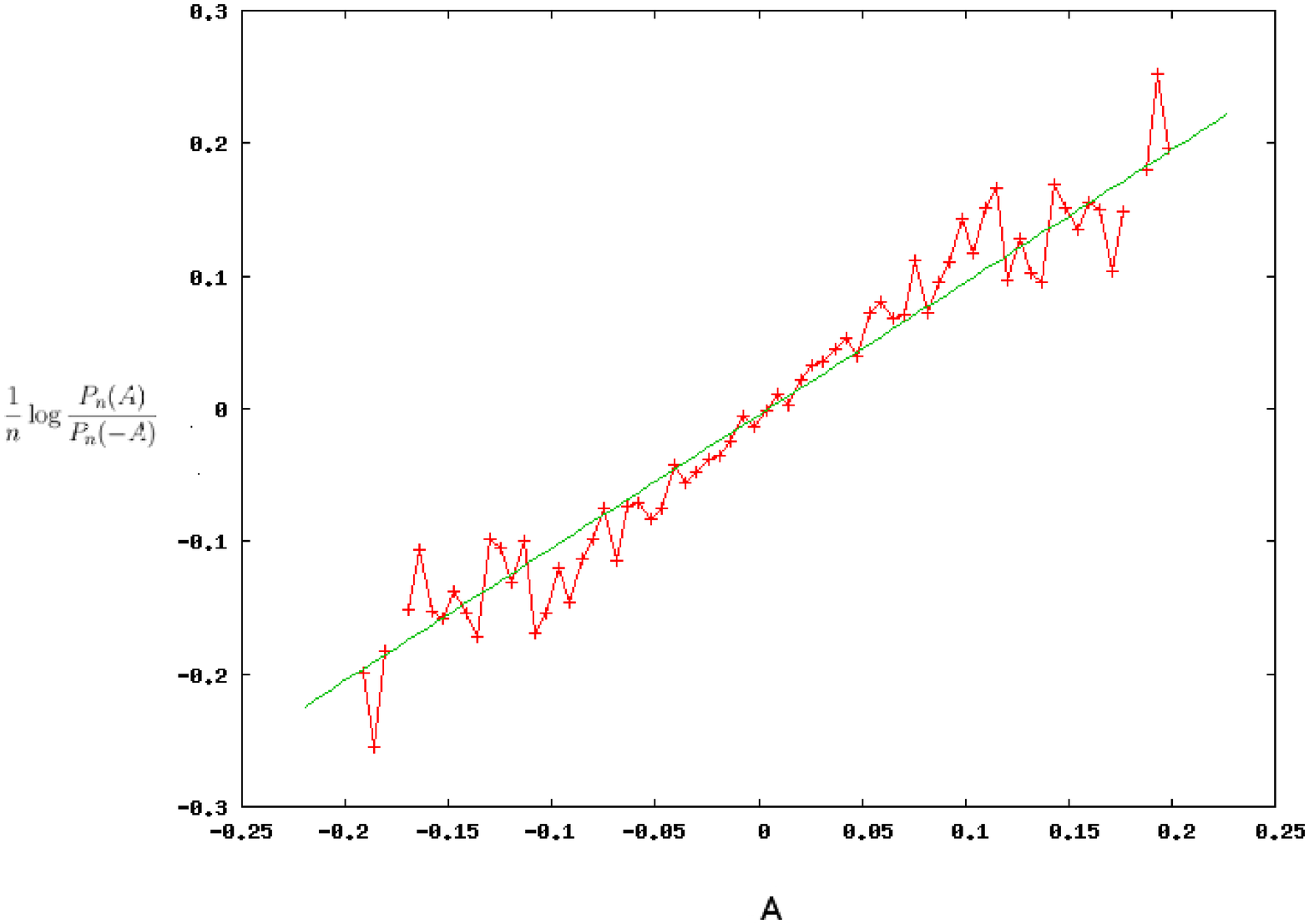}

\includegraphics[scale = 0.18]{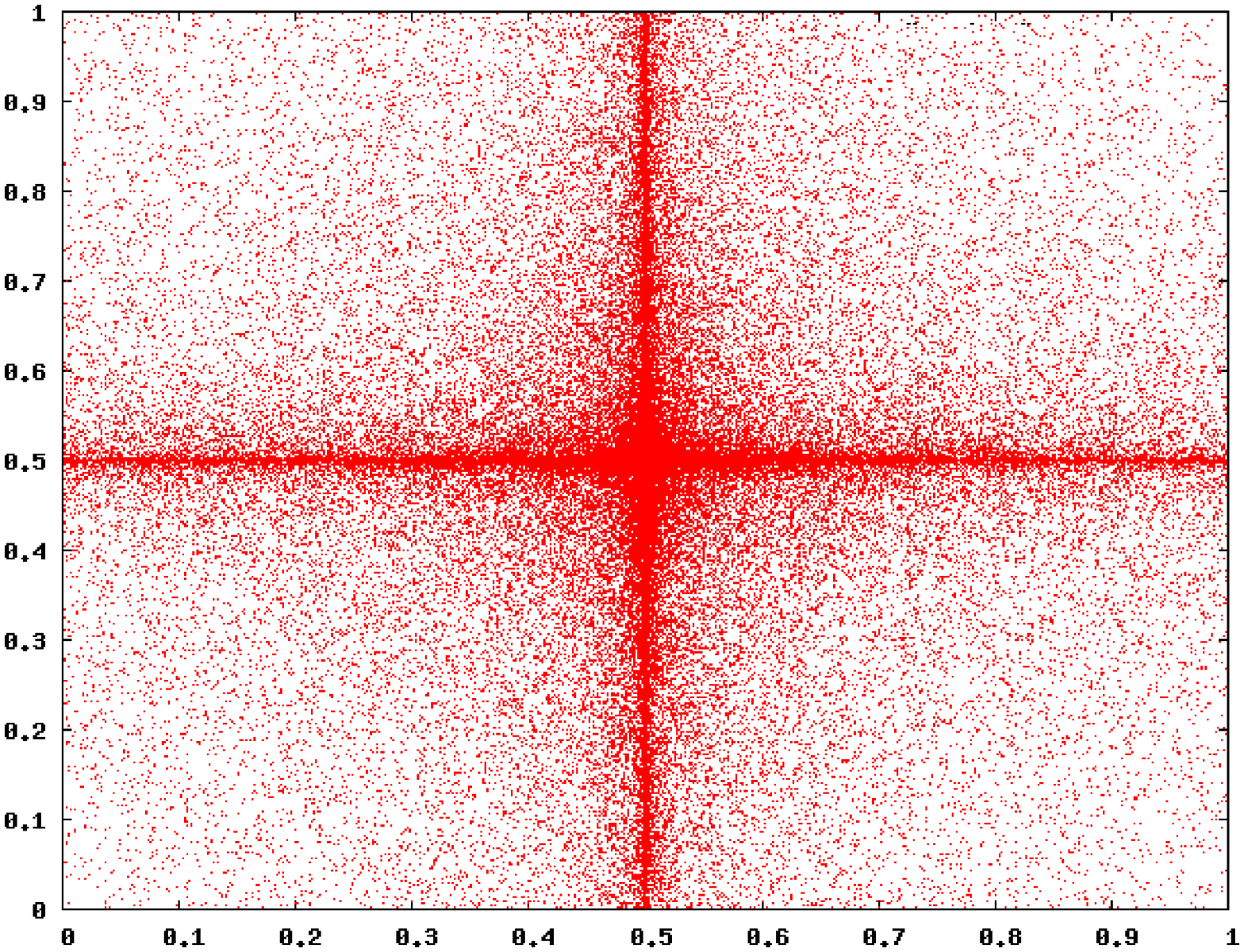}
\includegraphics[scale = 0.18]{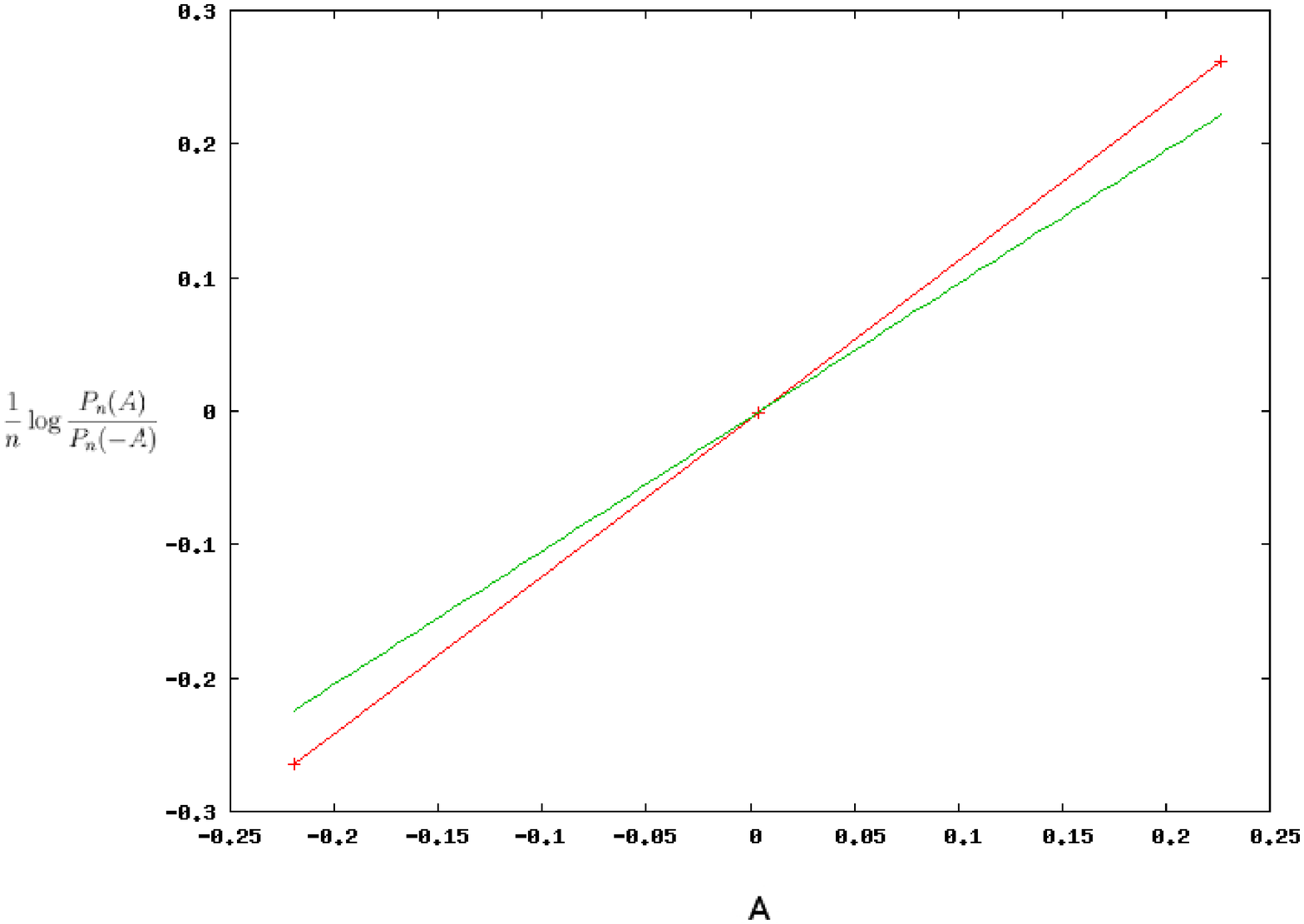}
\includegraphics[scale = 0.18]{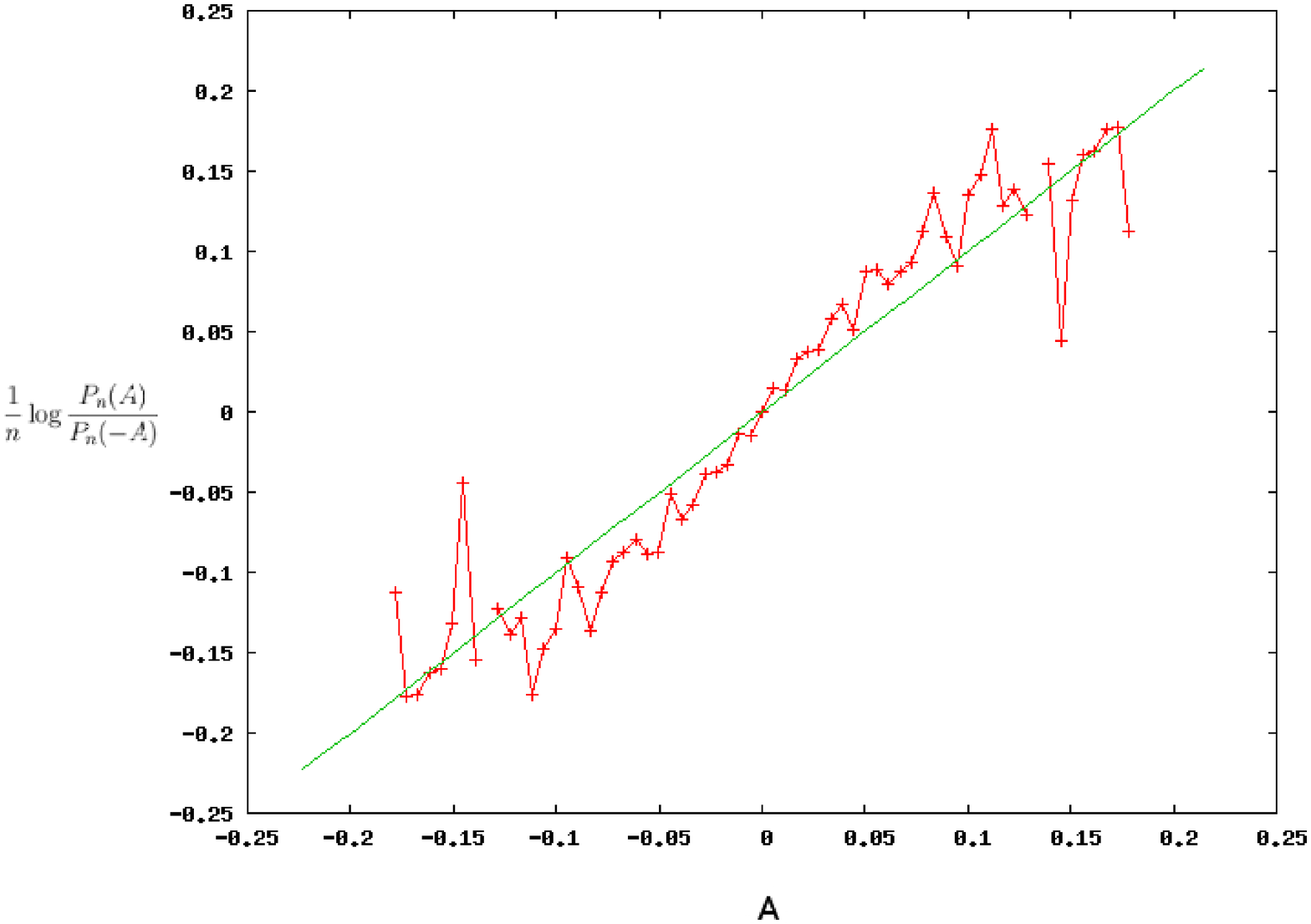}

\caption{ On the left, examples of different initial distributions. In the center, the FR is tested at step $n=1$ and $\ell=0.2$. Note that the FR is fulfilled just for the uniform distribution. On the right, the asymptotic FR gets restored for any initial condition, $n=40$ and $\ell=0.2$.  }\label{initdistr}
\end{figure}

We tested the robustness of our results by also varying the initial distribution, cf. Fig. \ref{initdistr}. 
Our analysis shows that the behaviour of the FR for small $n$ depends on the chosen initial distribution (for initial non-uniform distributions the FR is no longer automatically satisfied for $n=1$, as illustrated in Fig. \ref{initdistr}), as is the case with the transient FR discussed in Ref. \cite{ESR2}, whereas the validity of the asymptotic fluctuation relation is restored for large $n$, regardless of the initial (non singular) distribution.
In particular, for a given $\ell$, the number of steps needed to attain the proper asymptotic behaviour was found to depend only weakly on the initial distribution. 

\section{Conclusions}
The steady state dynamics of our map has no fluctuations. However, the map enjoys a weak form of time reversibility, which could lead in 
principle to the validity of the transient FR. As illustrated in \cite{ESR2} for reversible systems whose initial distribution may 
be generated by a single infinitely long equilibrium trajectory, the transient Fluctuation Relation 
(\ref{eqfrz}) is an identity valid for all $n$. Differently, our dynamical system violates 
the equilibrium ergodic condition. In fact, at equilibrium, all trajectories are period-4 cycles and the phase space is fragmented 
in various invariant sets of positive volume.
Our results thus emphasize the importance of the ergodicity of the equilibrium state.
Moreover our numerical simulations indicate that the validity of Eq.(\ref{eqfrz}) is restored in the long time limit.
This interesting property may not be immediately related to the standard transient FR, but it may stem from the peculiar features of the microscopic dynamics and seems to lead to a new kind of asymptotic, but not steady state, FR.

\section*{Acknowledgements}
Matteo Colangeli acknowledges financial support from the Brazilian agency Conselho Nacional de Desenvolvimento Cient\'{i}fico e Tecnol\'{o}gico (CNPq). 
%%%%%%%%%%%%%%%%%%%%%%%%%%%%%%%%%%%%%%%%%%%%%%%%%%%%%%%%%%%%%%%%%%%%%%%%%%%%%%%%%%%%%%%%%%%%%%%%%%%%%%%%%%%%%%%%%%%%%%%

\end{document}